\documentclass[11pt]{article}
\usepackage[english]{babel}
\usepackage{amsmath}
\usepackage{amsthm}
\usepackage{amsfonts}
\usepackage{amssymb}
\usepackage{hyperref}
\usepackage{enumitem}
\usepackage{titlesec}
\usepackage{secdot}
\sectiondot{subsection}
\sectiondot{subsubsection}
\usepackage{algorithm}
\usepackage{algpseudocode}
\usepackage{fancyhdr}
\usepackage{color}
\usepackage{tcolorbox}
\usepackage{subcaption}
\usepackage{titlesec}
\usepackage{relsize}
\usepackage{mathtools}

\usepackage{pgf,tikz}
\usepackage{graphicx} 

\usepackage{xcolor}

\theoremstyle{definition}
\newtheorem{threm}{Theorem}

\title{A Score Filter Enhanced Data Assimilation Framework for Data-Driven Dynamical Systems}

\author{Jingqiao Tang, Ryan Bausback, Feng Bao, Guannan Zhang, Phuoc-Toan Huynh}
\date{}

\begin{document}
\maketitle

\begin{abstract}
We introduce a score-filter-enhanced data assimilation framework designed to reduce predictive uncertainty in machine learning (ML) models for data-driven dynamical system forecasting. Machine learning serves as an efficient numerical model for predicting dynamical systems. However, even with sufficient data, model uncertainty remains and accumulates over time, causing the long-term performance of ML models to deteriorate. To overcome this difficulty, we integrate data assimilation techniques into the training process to iteratively refine the model predictions by incorporating observational information. Specifically, we apply the Ensemble Score Filter (EnSF), a generative AI-based training-free diffusion model approach, for solving the data assimilation problem in high-dimensional nonlinear complex systems. This leads to a hybrid data assimilation–training framework that combines ML with EnSF to improve long-term predictive performance. We Shall demonstrate that EnSF-enhanced ML can effectively reduce predictive uncertainty in ML-based Lorenz–96 system prediction and the Korteweg–De Vries (KdV) equation prediction. 
\end{abstract}

\section{Introduction}\label{sec1}

Machine learning (ML) has emerged as a powerful tool across a wide range of impactful domains, including computer vision~\cite{Chai2021, mao2024optimizing, Radosavovic2020}, natural language processing~\cite{Benes2017, Chai2019, Otter2021}, and weather forecasting~\cite{Allen2025, Beucler2021, Hakim2024}. More recently, neural network (NN)-based approaches have been developed as surrogate models for approximating and accelerating the simulation of dynamical systems governed by differential equations~\cite{Jiao2025, Lu2019, Lu2022}. Although a complete theoretical understanding remains an open question, these approaches can offer practical advantages over traditional solvers in some regimes. Firstly, once trained, they can provide rapid evaluations and can provide predictions at arbitrary spatial locations without extra interpolations. Secondly, for certain high-dimensional problems, learning-based surrogates can reduce the practical burden of spatial discretizations and thereby alleviate aspects of the curse of dimensionality.

However, despite numerous advantages of ML models, significant challenges remain. While such models often demonstrate strong short-term predictive accuracy, they frequently struggle to reliably capture long-term dynamics, particularly in chaotic systems~\cite{Pathak2018}. Without additional accurate inputs (e.g., new observational data for calibration) or continued retraining, these models can fail to generalize reliably to regimes far from the training distribution \cite{nguyen2004multiple}. For instance, the ``Aardvark Weather'' model in~\cite{Allen2025} demonstrated strong forecasting skill near the initial time but exhibited progressively degraded performance at longer lead times. Similarly, the DSDL framework for chaotic dynamics introduced in~\cite{Wang2024} accurately tracks the Lorenz–63 system over short horizons but shows substantial divergence from the true trajectory during extended time integration.

Another major challenge concerns the availability and quality of training data, which can significantly limit the development and performance of machine learning solvers. In practical applications, datasets are often limited in size, unevenly distributed, and sparsely sampled in space and time. Moreover, the available observations may only indirectly reflect the quantities of interest. In addition, training data are frequently affected by measurement noise and other sources of uncertainty. Such issues arise routinely in applications including weather and climate prediction, subsurface flow and transport, and ocean modeling, where data acquisition is costly and inherently uneven. As a consequence, even a well-trained surrogate model may inherit biases and systematic errors from the training data, ultimately leading to unreliable or unstable forecasts.

To mitigate error accumulation in long-term forecasts and reduce biases arising from imperfect training data, one common strategy is to collect additional measurements to recalibrate ML-based predictions. In practice, however, online measurements and observational data acquired during the forecasting phase are often partial, irregular, and contaminated by noise, making them insufficient to support reliable retraining. Moreover, repeatedly updating the model can be computationally expensive and may even degrade performance when the newly acquired data are limited or biased. These challenges motivate our alternative strategy that incorporates imperfect measurements directly into the prediction process. Instead of modifying the trained model, this approach updates the state estimate and its associated uncertainty online using observational information.
\vspace{0.5em}

In this paper, we introduce a ``data assimilation'' framework that integrates model predictions with observational data to reduce uncertainty and enhance the accuracy of machine learning (ML)-based forecasts. The underlying mathematical formulation of data assimilation is the optimal filtering problem, which seeks the best estimate of the target system state conditioned on the available observations. In particular, the optimal state estimate is characterized by the conditional expectation of the state variable given the observational information. The standard approach for solving the optimal filtering problem is the recursive Bayesian filter, which consists of the Kalman filter~\cite{LETKF_2020, evensen2009data,kalman1960} and the particle filter~\cite{Bao2019a, MTAC2012, pitt1999filtering, Andrieu2010, Gordon1993, Kang2018}. However, Kalman-type filters are often effective for high-dimensional problems, but they rely on linear or weakly nonlinear assumptions, which can limit their applicability in strongly nonlinear settings \cite{evensen2009data, LETKF_2007}. On the other hands, particle filters are, in principle, applicable to nonlinear and non-Gaussian problems, but they often suffer from severe weight degeneracy and the curse of dimensionality in high-dimensional state spaces because they are based on sequential Monte Carlo sampling~\cite{hu_vanLeeuwen_2021, Sny}.

To overcome the limitations of conventional methods, we adopt a score filter framework inspired by score-based diffusion models~\cite{song2019generative}, namely the Ensemble Score Filter (EnSF)~\cite{bao2024score, bao2024ensemble}. The EnSF aims to estimate the conditional probability density of the target model state, referred to here as the filtering density, and encodes its information in the corresponding score function. This score function then serves as a guiding term in a reverse-time diffusion process that steers the generative dynamics to produce samples consistent with the target distribution. Leveraging the strong expressive and generative capabilities of diffusion models \cite{YLiu2024, YSong2021}, extensive studies have demonstrated that EnSF can effectively track stochastic nonlinear dynamical systems in very high-dimensional spaces with high accuracy and computational efficiency \cite{bao2024ensemble, bao2025nonlinear, UF_2026, Bao2025a}. 

Building on its strong capability to filter uncertainties in predictive models, we integrate the EnSF framework with ML models to improve the predictive accuracy of dynamical systems learned from data. In this work, we consider two representative ML frameworks to demonstrate the performance of the EnSF-enhanced prediction: the long short-term memory (LSTM) network and operator-learning models.
The LSTM network, a gated recurrent neural network (RNN) designed for sequential and time-series modeling, naturally provides a step-by-step state transition map that propagates the system state forward in time. Its gating mechanism enables the network to capture long-term temporal dependencies while alleviating vanishing gradient issues, making it well suited for modeling dynamical systems.
In contrast, operator-learning models aim to learn mappings between function spaces and are typically trained to predict solutions over an entire spatiotemporal domain in a single forward pass. While effective for approximating complex dynamical systems, this global-in-time prediction structure does not naturally provide a recursive state-transition mechanism for time evolution. To match step-by-step state propagation, we adopt RNN-based architectures that advance the system state sequentially in time \cite{Goodfellow2016, Graves2012Supervised, lipton2015}.

The rest of the paper is organized as follows. In Section~\ref{Sec2:ML_Framework}, we review the machine learning components underlying our methodology, including the LSTM network and the Deep Operator Network (DeepONet) as a representative operator-learning framework. In Section~\ref{Sec3:DA_Framework}, we introduce the general data-assimilation setting, present the EnSF method, and discuss how it can be integrated with ML models to improve long-term predictive accuracy. In Section~\ref{Sec4:Results}, we present numerical experiments on LSTM-based Lorenz–96 prediction and DeepONet-based KdV equation solution prediction, and compare the EnSF approach with the ensemble Kalman filter under various uncertainty scenarios. Finally, we conclude the paper with a summary of the main findings.
 
\section{Overview of the Underlying Machine Learning Models}
\label{Sec2:ML_Framework}
In this section, we briefly review two machine learning frameworks adopted as surrogate models: the Long Short-Term Memory (LSTM) network and the operator-learning approach, in particular the Deep Operator Network (DeepONet). 

\subsection{Long-Short Term Memory}
\label{2.1:LSTM}

LSTMs are gated recurrent neural networks (RNN) introduced by the authors in~\cite{Hochreiter1997}. The key feature of RNN is that the current hidden state depends on the current input and the previous hidden state~\cite{Goodfellow2016}. However, unlike the vanilla RNN, the LSTM uses gated self-loops, which are known as cell states, to preserve information over long horizons and mitigate vanishing gradients~\cite{Goodfellow2016}. A compact form of the architecture of the LSTM is given by
\begin{equation}
\label{LSTM_solver}
\begin{array}{ll}
i_t = \sigma(W_i x_t + U_i h_{t-1} + b_i) & \text{(input gate)},\\
f_t= \sigma(W_f x_t + U_f h_{t-1} + b_f) & \text{(forget gate)},\\
\tilde{c}_t = \tanh(W_c x_t + U_c h_{t-1} + b_c)  & \text{(candidate cell)},\\
c_t = f_t \odot c_{t-1} + i_t \odot \tilde{c}_t & \text{(cell update)},\\
o_t = \sigma(W_o x_t + U_o h_{t-1} + b_o) & \text{(output gate)},\\
h_t = o_t \odot \tanh(c_t) & \text{(hidden state)},
    \end{array}
\end{equation}
where $x_t$ is the input at the time instance $t$, $h_{t-1}$ is the previous hidden state, $c_{t-1}$ is the previous cell, and $\odot$ is the element-wise product. $W_k, U_k, b_k$ are the learned parameters for $k=i, j, c, o$ and $\sigma$ is the logistic sigmoid. The forget gate $f_t \in (0,1 )$ controls how much of $c_{t-1}$ is carried forward, while the input and output gates regulate writing to and reading from the cell. We refer readers to~\cite{Gers2000,sherstinsky2020fundamentals} for a comprehensive overview of LSTMs.

\begin{figure}[h!]
\centering
\centerline{\includegraphics[width=0.8\linewidth,height=20pc]{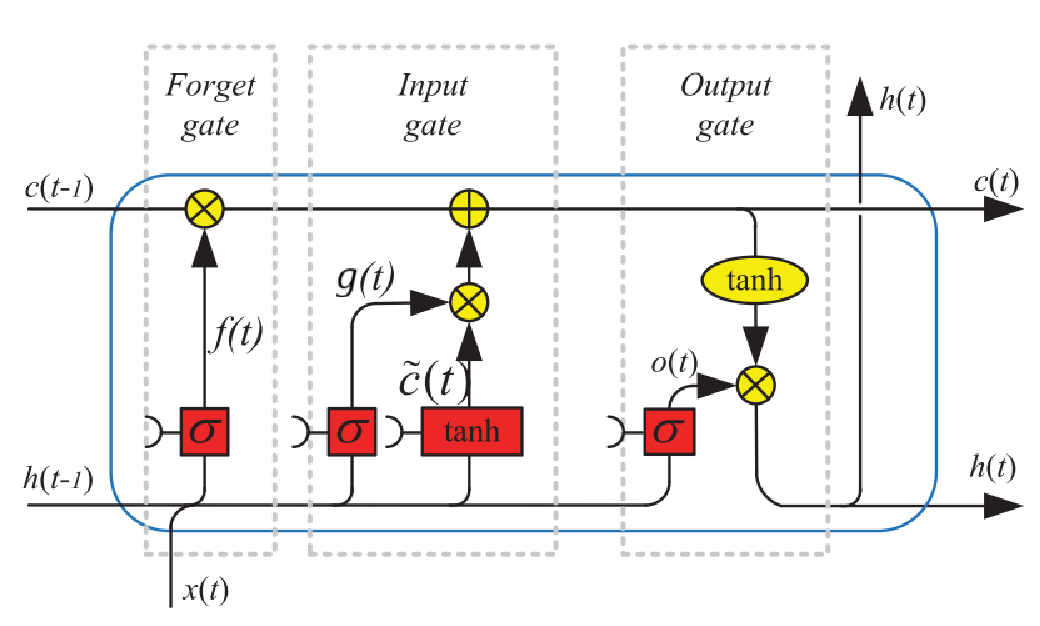}}
\caption{A block diagram of an LSTM cell, which takes recurrent input $h^{t-1}$, gated self-loop input $c^{t-1}$, and current state input vector $x^t$. These then flow through the forget gate $f^t$, the dual input gates $g^t$ and $\tilde{c}^t$, and output gate $o^t$ when making predictions for sequential and iterative data \cite{yu2019review}.}
\label{fig:lstm}
\vspace{-0.4cm}
\end{figure}

By adding this gated self-loop into the existing RNN structure, LSTM is able to more easily learn the long term dependencies within the sequential data and avoid losing critical information provided by small gradient values \cite{Goodfellow2016}. 

\subsection{Deep Operator Networks for Operator Learning}
\label{2.2:DeepONet}
Operator learning seeks to approximate an operator $G$ that maps an input function $u$ in a function space $V$ to an output function $G(u)$~\cite{Lu2019}. Unlike traditional supervised deep learning on finite-dimensional vectors, the input is an entire function, represented by its values 
$\{u(x_j)\}^m_{j=1}$ at $m$ sensor locations $\{x_j\}^m_{j=1}$, and the output is queried at point 
$y$ in the target domain \cite{Kovachki2024}. Given the discretized input $u(x)$ and a query point  $y$, the operator-learning architecture returns an approximation of $G(u)(y)$~ \cite{Lu2021}. Several neural operator frameworks have been proposed in this spirit; in this work, we focus on the DeepONet architecture and tailor it to the needs of data assimilation.

 DeepONet is a well-known method for solving the operator learning problem which was first developed in~\cite{Lu2019}. DeepONet takes $u$ and $y$ as inputs to separate, parallel neural networks, called the ``branch" and ``trunk" nets \cite{Lu2019}. The branch net takes the discretization of the function $u$, $\{u(x_0),..., u(x_m)\}$, evaluated at the $m$ sensor points $\{x_1, \hdots, x_m\}$, while the trunk net takes a single $y$ value in the operator output domain \cite{Lu2022}. These parallel networks then output vectors $[\beta_0\left(\pmb{u}\right),\hdots,\beta_p\left(\pmb{u}\right)]^T$ and $[\tau_0(y),\hdots,\tau_p(y)]^T$,  respectively. It is crucial that these output vectors have the same dimension $p$ so that they can be combined by the dot product to produce $G(u)(y)$. Then, the output of the DeepONet, the approximate operator of $G(u)$, is defined as:
\begin{equation}
\label{output_DeepONet}
\hat{G}(u)(y) = \sum_{k=1}^p\beta_k(\pmb{u})\tau_k(y).
\end{equation}

The theoretical foundation for the DeepONet was first established by the authors in~\cite{Chen1995}: 

\begin{threm}[\textbf{Universal Approximation Theorem for Operators}~\cite{Chen1995}]
\label{univapproxop}
Suppose $\sigma$ is a continuous non-polynomial function, $X$ is a Banach space, $K_1 \subset X$, $K_2 \subset \mathbb{R}^d$ are compact sets, $V$ is a compact subset of $C(K_1)$, and $G$ is an operator mapping $V$ to $C(K_2)$. Then for $\epsilon > 0$, there are positive integers $n$, $p$, $m$, constants $c^k_i$, $\xi^k_{ij}$,$\theta^k_i$ $\in \mathbb{R}$, $w_k \in \mathbb{R}^d$, $x_j \in K_1$, $i = 1, \hdots, n$,
$k = 1,\hdots, p$ and $j = 1,\hdots, m$, such that
$$\left|G(u)(y) - \sum^p_{k=1} \underbrace{\sum^n_{i=1}c^k_i \sigma\left(\sum^m_{j=1} \xi_{ij}^ku(x_j)+\theta_i^k\right)}_{\text{Branch}}\underbrace{\sigma(w_ky+ \zeta_k)}_{\text{Trunk}}\right| < \epsilon, $$
holds for all $u \in V$ and $y \in K_2$. Here, $C(K)$ is the Banach space of all continuous
functions defined on $K$ with norm $\vert\vert f \vert\vert_{C(K)} = \max\limits_{x \in K} \vert f(x) \vert.$
\end{threm}

This theorem proves that a shallow branch–trunk pair can approximate any continuous operator. However, it also restricts both sub-networks to have a single hidden layer and identical structure. 
Subsequently, the authors in~\cite{Lu2021} removed these constraints and proved an analogous approximation guarantee for arbitrary (potentially deep and asymmetric) branch and trunk architectures. 
We summarize their extension as Theorem~\ref{univapproxop2} below.

\begin{threm}[\textbf{Generalized Universal Approximation Theorem for Operators},~\cite{Lu2021}]
\label{univapproxop2}
Suppose $X$ is a Banach space, $K_1 \subset X$, $K_2 \subset \mathbb{R}^d$ are two compact sets, respectively, $V$ is compact in $C(K_1)$, and $G: V \rightarrow C(K_2)$ is a nonlinear continuous operator. Then, for any $\epsilon >0$, there exist positive integers $m, p,$ continuous vector functions $\pmb{g}: \mathbb{R}^m \to \mathbb{R}^p$, $\pmb{f}: \mathbb{R} \to \mathbb{R}^p$, and $x_1, \hdots, x_m \in K_1$, such that

$$\left|G(u)(y) - \langle \underbrace{\pmb{g}(u(x_1), u(x_2), ... u(x_m))}_{\text{Branch}}, \underbrace{\pmb{f}(y)}_{\text{Trunk}} \rangle\right| < \epsilon, $$
holds for all $u \in V$ and $y \in K_2$, where $\langle \cdot, \cdot \rangle$ denotes the dot product in $\mathbb{R}^p$. Furthermore, the functions $\pmb{g}$ and $\pmb{f}$ can be chosen as diverse classes of neural networks, which satisfy the classical universal approximation theorem of functions, for examples, (stacked/unstacked) fully connected neural networks, residual neural networks and convolutional neural networks.
\end{threm}

\subsection{Recurrent Deep Operator Network}
\label{2.3:Combination}
To learn time-dependent dynamical systems, we integrate an LSTM structure into DeepONet to obtain an alternative solver formulation, which we refer to as the Recurrent DeepONet (R-DeepONet) in this paper and which is also known as “S-DeepONet” in \cite{S-DeepONet}. This is justified by the universal approximation property of Theorem~\ref{univapproxop2}. Assume that $v_t(x)$ denotes the solution of a given PDE at time $t$ and spatial position $x$. In the standard setting, DeepONet is used as a one-step solver: the solution at time $t_{n+1}$ is computed from the solution at the previous time step $t_n$. To incorporate information from multiple past states, we replace the branch network by an LSTM that processes several previous time steps as input. Specifically, we define the R-DeepONet solution operator as
\begin{equation}
\label{DeepONet_Solver}
    G_{\text{R-DeepONet}} : U_n \mapsto v_{t_{n+1}},  
\end{equation}
where
\begin{align*}
U_n := \bigl(v_{t_{n-k}}, \dots, v_{t_n}\bigr), \; k\geq 0.
\end{align*}
Since we use the same spatial discretization across the entire trajectory, the trunk input 
$y$ is taken to be a sensor location $x$ at which each snapshot $v_{t_n}(x)$ is evaluated. Accordingly, the R-DeepONet outputs a one-step-ahead prediction, 
\begin{align*}
v_{t_{n+1}}(y) = G_{\text{R-DeepONet}}(U)(y),
\end{align*}
i.e., a single point on the next cross-section of the operator trajectory. Therefore, in this recurrent setting the sensor locations $x$ and the output evaluation points $y$ coincide. The loss function remains the same as in the vanilla DeepONet, since the output quantity being fitted is unchanged; only the input is reinterpreted to support step-by-step solution tracking.

The LSTM solver~\eqref{LSTM_solver} and the R-DeepONet solution operator~\eqref{DeepONet_Solver} can be used as alternatives to traditional numerical solvers, since they are mesh-free and require little manual adjustment when the computational domain is irregular or contains many sharp corners. 

\begin{figure}[h!]
\centerline{\includegraphics[width=0.7\linewidth, height=7pc]{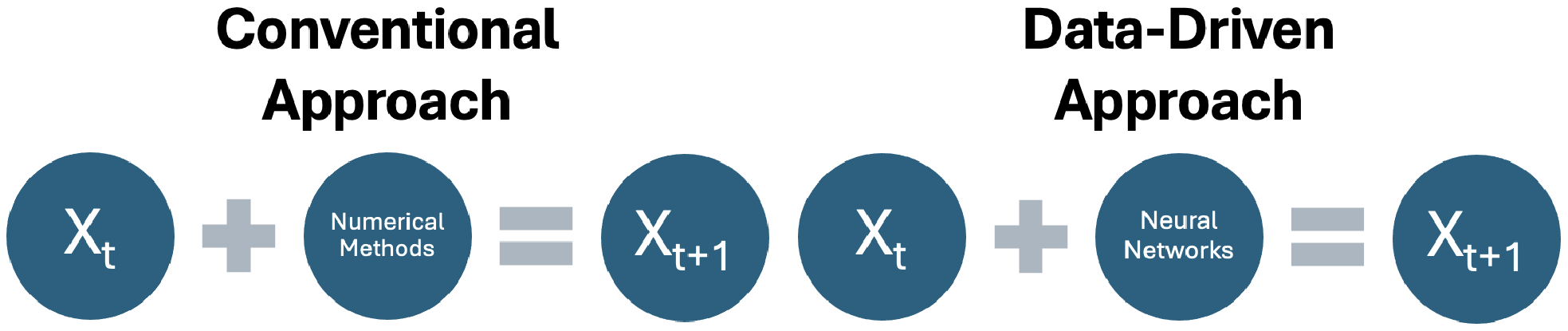}}
\caption{The traditional approach of numerically solving the governing equations versus data-driven approach on predicting the spatio-temporal models. }
\label{fig:system_tracking}
\vspace{-0.2cm}
\end{figure}

However, because these machine-learning solvers are constructed purely from data, their predictions can become unreliable when the available data do not accurately reflect the behavior of the true state. This mismatch is common in real-world settings where observations are noisy, sparse, and sometimes strongly nonlinear. Moreover, when these models are used to predict the evolution over a very long time window, small forecast errors can accumulate and lead to substantial deviations. Even with abundant training data, residual model error can persist, and systematic improvements are often difficult because, unlike traditional numerical schemes, these solvers do not provide transparent analytic structure. These limitations motivate the use of data assimilation, which fuses model forecasts with observational data to sequentially correct the state and reduce uncertainty.

In the next section, we present the general framework of the data assimilation problem and recall two filtering methods that will be employed in our study: the classical Ensemble Kalman Filter (EnKF) approach  and the recently developed  Ensemble Score Filter (EnSF) method. These methods are then combined with the LSTM-based networks to provide a hybrid data assimilation–driven prediction framework.
\section{Data Assimilation Framework and the Ensemble Score Filter}
\label{Sec3:DA_Framework}
We consider the following stochastic dynamical system given by
\vspace{-0.3cm}
\begin{equation} \label{eqn:da_forward_equation}
X_{t} = f(X_{t-1}, \omega_{t-1}), 
\end{equation}
where $f:\mathbb{R}^d \times \mathbb{R}^k \mapsto \mathbb{R}^d$ is the (possibly nonlinear) dynamical model, \; $X_t \in \mathbb{R}^d$ is the current state of the dynamical system at time instant $t$ defined by $f$, and $\omega_t \in \mathbb{R}^k$ is a random variable representing the model uncertainty. This uncertainty can arise from many sources, including noisy or inaccurate data, unknown physical parameters, and modeling errors due to incomplete knowledge of the underlying physics. As the forward model propagates in time using the solver $f$, these uncertainties and errors naturally accumulate, leading to progressively less accurate predictions. To mitigate this deterioration in accuracy, we incorporate observational data to correct the model predictions and improve the state estimates. Such data are often collected through the following observation process:
\vspace{-0.2cm}
\begin{equation} \label{eqn:da_observation_equation}
Y_t = h(X_t)+E_t,
\vspace{-0.1cm}
\end{equation}
where $Y_t \in \mathbb{R}^r$ is the observation, and $E_t \sim N(0, \Sigma_t)$ is a Gaussian noise with covariance matrix $\Sigma_t$.

The goal of data assimilation is to give the best possible estimate of the state $X_t$ given the observational data up $Y_{1:t}$. The mathematical framework used to solve the data assimilation problem is the \textit{optimal filter}, which seeks the best estimation in the form of the conditional expectation of $X_t$ given $Y_{1:t}$. Specifically, the optimal filter $\hat{X}_t$ of $X_t$ is given by:
\begin{align*}
    \hat{X_t} := \mathbb{E}[X_t | \mathcal{Y}_t],
\end{align*}
where $\mathcal{Y}_t := \sigma{(Y_{1:t})}$ is the 
$\sigma$-algebra generated by the observations $\{Y_1, \hdots, Y_t\}$. Alternatively, one can approximate such conditional expectation by computing the conditional probability density function $p(X_t | Y_t)$. This is done through the Bayesian filtering approach which consists of two steps: the prediction step and the update step. 

In the prediction step, we use the Chapman-Kolmogorov formula to compute the prior filtering density of the state variable at time $t$:
$$p(X_{t} | Y_{t-1}) = \int p(X_t | X_{t-1})p(X_{t-1} | {Y}_{t-1}) dX_{t-1},$$
where $p(X_t | X_{t-1})$ is the transition probability derived from the state dynamics~\eqref{eqn:da_forward_equation}. 

In the update step, the posterior  filtering density is computed by using the following Bayesian inference scheme:
$$p(X_{t} | {Y}_{t}) \propto p(X_{t}|{Y}_{t-1}) p(Y_t | X_t),$$
 where $p(Y_t | X_t)$ is the likelihood function that quantifies the probability of the observed data given the current state:
$$p(Y_t | X_t) \propto \exp\left[-\frac{1}{2}(h(X_t)-Y_t)^T \Sigma_t(h(X_t)-Y_t)\right].$$ 

The prediction step and the update step form a recursive model-data integration framework that dynamically improves the accuracy of the system's state. However, for high-dimensional problems it is challenging to accurately represent the filtering distribution and to effectively assimilate new observational data. In what follows, we introduce an effective ensemble score filter (EnSF) approach for solving the data assimilation problem.

\subsection{Ensemble Score Filter}
\label{3.2:EnSF}

The Ensemble Score Filter (EnSF) is a novel approach for solving high-dimensional nonlinear filtering problems, building on the score-based filtering framework developed in~\cite{bao2024score}. In contrast to particle filters and EnKF, which represent the filtering distribution with a finite ensemble of particles, EnSF encodes the recursively updated filtering density in its score function \cite{bao2024score}. This method is inspired by score-based diffusion models, a class of generative machine learning methods that represent a target distribution through its score and generate state samples by simulating a reverse-time stochastic differential equation (SDE)~\cite{sohl2015deep}.

In the score-based diffusion model framework, two stochastic processes are defined as follows:
\begin{equation}\label{Diffusion_Model}
\fontsize{11pt}{11pt}\selectfont
\begin{array}{ll}
    d Z_{\tau} = &  b_{\tau} Z_{\tau} d_{\tau} + \sigma_{\tau} dW_{\tau}, Z_0 = Q_0(X), \;\; \text{(Forward SDE)}\\
    d Z{\tau} = & \big( b_{\tau} Z_{\tau} - \sigma_{\tau}^2 S(Z_{\tau}, \tau) \big)d_{\tau} + \sigma_{\tau} d\overleftarrow{W}_{\tau}, Z_t = N(0, \pmb{I}_d). \;\;\text{(Reverse SDE)}
\end{array}
\end{equation}
The first equation is a forward stochastic differential equation (SDE) that propagates a given target data distribution $Z_0= Q_0(X)$, such as the filtering distribution in data assimilation, toward the standard Gaussian distribution over the pseudo time interval $[0, 1]$. The second equation is the corresponding time-reversed SDE, which transports the standard Gaussian random variable $Z_{1} = N(0, \pmb{I}_d)$ back to the original data distribution $Z_0$, where $\int \cdot d\overleftarrow{W}_{\tau}$ is a backward stochastic integral \cite{Bao_AA20, BaoC20142, BaoF20182}. The drift coefficient $b$ and diffusion coefficient $\sigma$ in the forward SDE are pre-defined deterministic functions. In this work, we let $b_{\tau} = \dfrac{d \log \alpha_{\tau}}{ d \tau }$ and $\sigma_{\tau} = \sqrt{\dfrac{d \beta_{\tau}^2}{d \tau} - 2 \dfrac{d \log \alpha_{\tau} }{d \tau} \beta_{\tau}^2 }$ with $\alpha_{\tau} = 1 - \tau$ and $\beta_{\tau} = \tau$ (see \cite{bao2024ensemble}).

The key component that enables the reverse SDE in Eq.~\eqref{Diffusion_Model} to transport samples from the standard Gaussian distribution to the target data distribution is the \textit{score function} $S(Z_{\tau}, \tau)$, defined by 
\begin{equation}\label{Def:Score}
S(Z_{\tau}, \tau) = \nabla_z \log Q_{\tau}(Z_{\tau}),
\end{equation}
where $Q_{\tau}$ is the distribution of the forward SDE solution $Z_{\tau}$. In other words, the score function $S$ acts as a guide that steers the reverse-time sampling process in the diffusion model framework back toward the original data distribution. Thus, if we have access to the score $S(Z_\tau, \tau)$ associated with $P(X_{t}|Y_{1:t})$, we can simulate the reverse-time SDE in~\eqref{Diffusion_Model} to generate samples following such filtering distribution. 

To complete the data assimilation task, we need to construct the scores $S_{t|t-1}$ and $S_{t|t}$ corresponding to the prior filtering density $P(X_{t}|Y_{1:t-1})$ and the posterior filtering density $P(X_t|Y_{1:t})$, respectively. To obtain the score $S_{t|t-1}$, we use the definition of marginal and conditional probability density functions to write the target density $Q_\tau(Z_\tau)$ as
$$Q_\tau(Z_\tau) = \int_{\mathbb{R}^d}Q_{\tau | 0}(Z_\tau|Z_0)Q_0(Z_0)dZ_0,$$
and thus
\begin{equation}
\fontsize{11pt}{11pt}\selectfont
\begin{array}{l}
S_{t|t-1}(Z_\tau, \tau) = \nabla_{Z_\tau} \log \left(\mathlarger{\int}_{\mathbb{R}^d}Q_{\tau | 0}(Z_\tau|Z_0)Q_0(Z_0)dZ_0\right) \vspace{0.2cm} \\
=\dfrac{1}{\mathlarger{\int}_{\mathbb{R}^d}Q_{\tau | 0}(Z_\tau|Z'_0)Q_0(Z'_0)dZ'_0}\mathlarger{\int}_{\mathbb{R}^d}-\frac{Z_\tau -\alpha_\tau Z_0 }{\beta^2_\tau}Q_{\tau | 0}(Z_\tau|Z_0)Q_0(Z_0)dZ_0 \vspace{0.1cm} \\
= \mathlarger{\int}_{\mathbb{R}^d}-\frac{Z_\tau -\alpha_\tau Z_0 }{\beta^2_\tau}w(Z_\tau,Z_0)Q_0(Z_0)dZ_0,
\end{array}
\end{equation}
where 
$$w(Z_\tau,Z_0) := \frac{Q_{\tau | 0}(Z_\tau|Z_0)}{\int_{\mathbb{R}^d}Q_{\tau | 0}(Z_\tau|Z'_0)Q_0(Z'_0)dZ'_0}.$$
After the score $S_{t|t-1}$ is available, we compute the posterior score using the following relation:
\begin{equation}
\label{update_postscore}
S_{t|t}(Z_\tau, \tau; \theta)= S_{t|t-1}(Z_\tau, \tau; \theta) + v(\tau)\nabla_z\log P(Y_{t}|Z_{\tau}),
\end{equation}
 where $P(Y_{t}|Z_{\tau})$ is analytically defined as
 \begin{equation}
 p(Y_t | Z_t) \propto \exp\left[-\frac{1}{2}(h(Z_\tau)-Y_t)^T \Sigma (h(Z_\tau)-Y_t)\right],
 \end{equation}
and $v(\tau)$ is a monotonically decreasing damping function with $v(0)=1$ and $v(1)=0$ \cite{bao2024score}. The damping function acts like a step size in stochastic gradient descent and in practice, we use $v(\tau) = 1 -\tau$. This results in the likelihood having very little impact when $\tau \approx 1$, but as $\tau$ decreases, the diffusion term becomes less dominant, and the likelihood is gradually injected via the drift term \cite{bao2024ensemble}.

In the first work, the score functions are approximated by using neural network models~\cite{bao2024score}. However, this approach is expensive and inefficient as the models need to be retrained every time we receive new observational data. To overcome this challenge, the authors in~\cite{bao2024ensemble} propose the EnSF which directly discretizes the scores $S_{t|t-1}$ and $S_{t|t}$. Specifically, we first draw $J$ samples $\{x^j_{t-1| t-1}\}_{j=1}^J$ from $P(X_{t-1}|Y_{1:t-1})$ by solving the reverse-time SDE~\eqref{Diffusion_Model} using $S_{t-1|t-1}$. We then simulate the sample $\{x^j_{t|t-1}\}^J_{j=1}$ using the forward solver $f$ in~\eqref{eqn:da_forward_equation} and extract from it a mini-batch $\{x^{j_n}_{t|t-1}\}_{n=1}^N$ to provide an approximation for the score $S_{t|t-1}$ \cite{EnSF_Scalable_2024, IEnSF_2025}: 
\begin{equation}
\fontsize{11pt}{11pt}\selectfont
\label{estimate_priorscore}
S_{t|t-1}(Z_\tau, \tau) \approx \hat{S}_{t|t-1}(Z_\tau, \tau): = \sum_{n=1}^N-\frac{Z_\tau -\hat{\alpha}_\tau x^{j_n}_{t|t-1}}{\hat{\beta}^2_\tau}\bar{w}(Z_\tau,x^{j_n}_{t|t-1}),
\end{equation}
 where $\bar{w}$ is an approximation of the weight function $w$ given by:
 \begin{align*}
 \bar{w}(z, x^{j_n}_{t\vert t-1}) := \frac{Q_{t | 0}(z|x^{j_n}_{t|t-1})}{\sum\limits_{m=1}^NQ_{t | 0}(z|x^{j_m}_{t|t-1})}.
 \end{align*}
 Finally, we use the formula~\eqref{update_postscore} with the estimated score $\hat{S}_{t|t-1}$ in~\eqref{estimate_priorscore} to obtain an approximation of the posterior score:
 \begin{equation}
 \fontsize{11pt}{11pt}\selectfont
\label{estimate_postscore}
S_{t|t}(Z_\tau, \tau; \theta) \approx \hat{S}_{t|t}(Z_{\tau}, \tau; \theta) = \hat{S}_{t|t-1}(Z_\tau, \tau; \theta) + v(\tau)\nabla_z\log P(Y_{t}|Z_{\tau}).
\end{equation}
With the approximate posterior score $\hat{S}_{t|t}$ in~\eqref{estimate_postscore}, we can generate a set of samples $\{x^j_{t|t}\}^J_{j=1}$ that follow
the posterior filtering density $p\left(X_t|Y_{1:t}\right)$ from the Gaussian distribution through the reverse-time SDE~\eqref{Diffusion_Model}. The posterior state at time $t$ is then given by
\begin{equation}
    \hat{X}_t = \dfrac{1}{J}\sum\limits^J_{j=1}x^{j}_{t|t}.
\end{equation}

\subsection{Integrating Data Assimilation with Machine Learning Models}
\label{3.3:DA+ML}
In the numerical experiments, we consider two recursive prediction modes based on a $K$-step input window. Let $f(\cdot)$ denote the trained LSTM-based model (e.g., LSTM or R-DeepONet). For each $k\ge K$, we define:
\begin{equation}\label{eqn:LSTM_SSP}
(SSP)_k = f\!\left(\{Y_i\}_{i=k-K}^{k-1}\right),
\end{equation}
where $\{Y_i\}_{i=k-K}^{k-1}$ are the \emph{true} previous states, and
\begin{equation}\label{eqn:LSTM_LTP}
(LTP)_k = f\!\left(\{(LTP)_i\}_{i=k-K}^{k-1}\right),
\end{equation}
where $(LTP)_i = Y_i$ for $0\le i\le K-1$, and $(LTP)_k$ is generated recursively for $k\ge K$.
We refer to the first approach as \emph{single-step prediction (SSP)} and the second as \emph{long-term prediction (LTP)}. Since LTP uses previously predicted states as inputs, model errors can accumulate over time, which motivates integrating data assimilation to reduce this drift and improve forecasts. The pseudo-code for combining data assimilation with the LTP-based solver is given in Algorithm~\ref{DA_LSTM_LTP}.
\begin{algorithm}[h!]
\caption{Data assimilation with an LTP-based forecast model}\label{DA_LSTM_LTP}
\begin{algorithmic}[1]
\State \textbf{Input:} training dataset; observation operator $g(\cdot)$; time window $K$;
ensemble size $J$; prior $P(Y_{t_0})$;  observation noise covariance $R$.
\State Train the ML forecast model (either LSTM or R-DeepONet) $f(\cdot)$ using the training dataset.
\State Sample initial ensemble $\{y^{(j)}_{0}\}_{j=1}^{J}\sim P(Y_{t_0})$.
\State Initialize $\{y^{(j)}_{k}\}_{j=1}^{J}$ for $k=1,\dots,K-1$ using a high-fidelity numerical solver.
\For{$n=K, K+1, \dots$}
\State \textbf{Forecast (LTP):} $\tilde{y}^{(j)}_{n} = f\!\left(\{y^{(j)}_{i}\}_{i=n-K}^{n-1}\right)$ for $j=1,\dots,J$.
\State \textbf{Observation model:} given $y_n$, compute $z_n = g(y_n)+\eta_n$, $\eta_n\sim\mathcal N(0,R)$.
\State \textbf{Analysis:} update $\{\tilde{y}^{(j)}_{n}\}_{j=1}^{J}$ to obtain analysis ensemble $\{\hat{y}^{(j)}_{n}\}_{j=1}^{J}$.
\If{Using EnKF}
    \State Compute predicted observations $\tilde z_n^{(j)} = g(\tilde y_n^{(j)})$ and apply EnKF update to form $\{\hat{y}^{(j)}_{n}\}_{j=1}^{J}$.
\ElsIf{Using EnSF}
    \State Estimate prior score $S_{n|n-1}$ from $\{\tilde{y}^{(j)}_{n}\}_{j=1}^{J}$ using~\eqref{estimate_priorscore}.
    \State Compute posterior score $S_{n|n}$ by combining $S_{n|n-1}$ with the likelihood using~\eqref{estimate_postscore}.
    \State Solve the reverse-time SDE in~\eqref{Diffusion_Model} to generate analysis ensemble $\{\hat{y}^{(j)}_{n}\}_{j=1}^{J}$.
\EndIf
\State State estimate at time step $t_n$: $\bar y_n = \frac{1}{J}\sum\limits_{j=1}^{J}\hat{y}^{(j)}_{n}$.
\State Set $y^{(j)}_{n} \leftarrow \hat{y}^{(j)}_{n}$ for the next filtering step.
\EndFor
\end{algorithmic}
\end{algorithm}


\section{Numerical Experiments}
In this section, we present a series of numerical experiments to demonstrate the effectiveness of the diffusion model–based ensemble score filter in reducing uncertainty in AI-based forward models and enhancing predictive accuracy. We first consider the Lorenz–96 system in Section \ref{4.1:Lorenz96}, a canonical benchmark for chaotic dynamics, using a standard LSTM model as the predictive solver. We then study the Korteweg–De Vries equation in Section \ref{4.2:KdV}, where an R-DeepONet surrogate is employed as an AI-based forward model for learning PDE dynamics. In both settings, we show that the proposed approach substantially improves long-horizon predictive performance and effectively reduces uncertainty in machine learning solvers. Moreover, our experiments demonstrate that the ensemble score filter can calibrate imperfectly trained surrogate models, leading to more reliable predictions even in data-limited or insufficiently trained regimes.

\label{Sec4:Results}
\subsection{Lorenz-96 Model}
\label{4.1:Lorenz96}
We consider the following (20-dimensional) Lorenz-96 model :
\begin{equation} \label{eqn:Lorenz96_Definition}
    \frac{dx_t^{i}}{dt} = (x_t^{i + 1} - x_t^{i - 2})x_t^{i - 1} - x_t^{i} + F, \; i=1, \hdots, 20,
\end{equation}
where $X_t = \left[x^{1}_t, \hdots, x^{20}_t\right]^T$ is the state variable. 
System~\eqref{eqn:Lorenz96_Definition} is closed by enforcing periodic boundary conditions with the dimension period, i.e.,
$x^{0} = x^{20}$, $x^{-1} = x^{19}$, and $x^{1} = x^{21}$ (equivalently, $x^{i+20}=x^{i}$ for all integers $i$).
The external forcing parameter is fixed at $F = 8$. For a random initial state $X_{0}$, we compute the numerical trajectory using a Runge–Kutta ODE solver from SciPy evaluated on a uniform temporal grid of $\Delta{t} = 10/T$, where $T=2500$, over the time interval $[0, 10]$. 

\subsubsection{Training configuration for the LSTM model} \label{sec:LSTM_model_training}
\begin{figure}[h!]
\centering
\includegraphics[width = 0.55\linewidth]{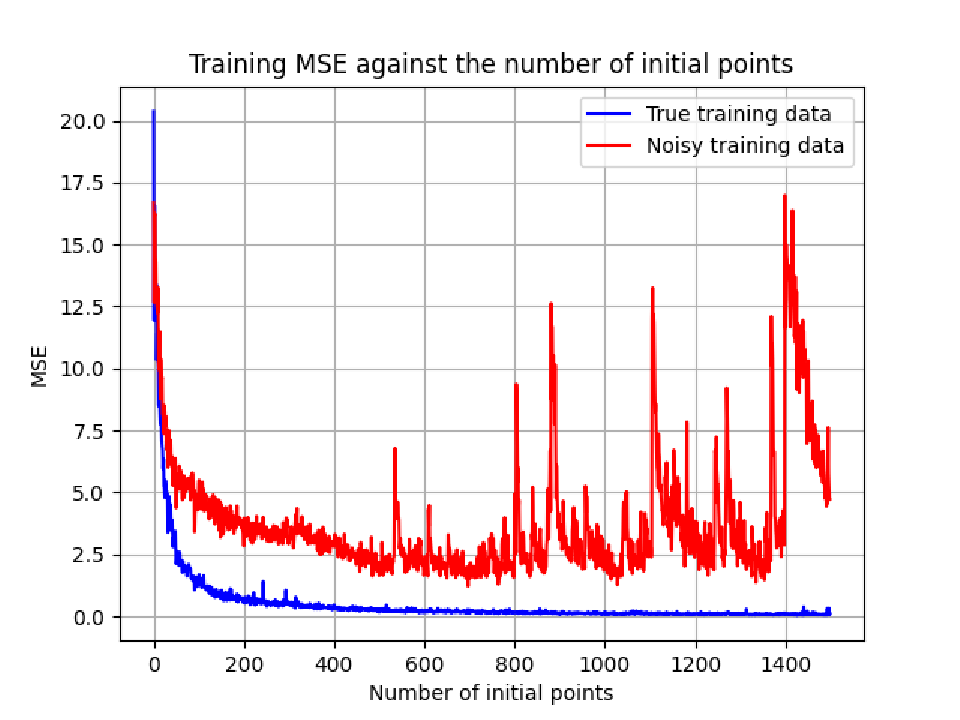}
\caption{LSTM training loss for clean and noisy datasets. The model trained on clean data achieves low training error with stable convergence, while training on noisy data leads to less stable convergence and higher residual error.}\label{fig:LSTM_training_MSE_total}
\label{sec:SSP_LTP}
\vspace{-0.2cm}
\end{figure}

The LSTM model used in this numerical experiment consists of five recurrent layers, each with 128 hidden units, and is trained using the Adam optimizer with a learning rate of $\eta = 10^{-3}$.A memory length of $20$ is employed, meaning that states from the previous $20$ time steps are used to predict the subsequent state.  For each training realization, we generate a random initial state $X_0$ and integrate the system to obtain a trajectory $\{X_t\}_{t=0}^{2500}$ using the numerical solver described above. From this trajectory, we extract supervised input--output pairs$\big(\{X_{t-19},\ldots,Y_t\},\, Y_{t+1}\big)$, $t=19,\ldots,2499$, so that the first target is $X_{20}$ and the last target is $X_{2500}$. To ensure sufficient training variability, this procedure is repeated $1500$ times with independently generated random initial states and their associated trajectories. The LSTM forward solver used in this section is trained on the resulting dataset for $200$ epochs. The mean squared error (MSE) training loss is reported in Figure~\ref{fig:LSTM_training_MSE_total}, where the blue curve illustrates the convergence behavior in the clean-data setting. 

In this work, we also consider a setting in which the training data are contaminated by noise, representing observational errors or potential adversarial perturbations. Specifically, the same LSTM model is trained using trajectories corrupted by additive Gaussian noise:
$
\tilde{X}_t = X_t + \xi_t, \; \xi_t\sim\mathcal{N}(0,0.2^2),
$
where the noise is applied independently to each component. The corresponding MSE loss curve is also shown in Figure~\ref{fig:LSTM_training_MSE_total} (red curve). Compared with the clean-data case, training with noisy data results in a higher loss magnitude and more pronounced oscillations, suggesting that noise can degrade LSTM training process and leads to increased uncertainty in predictive accuracy.

\begin{figure}[htb]
\centering
\includegraphics[scale = 0.5]{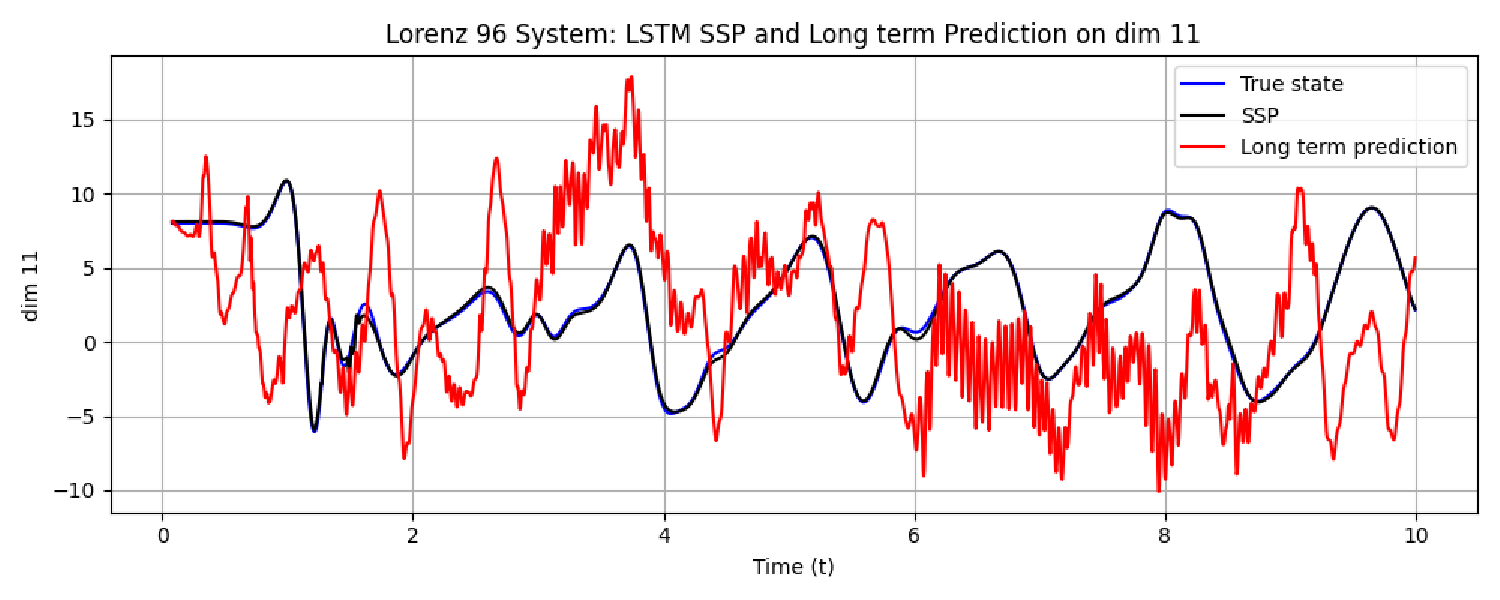}
\vspace{-1em}
\caption{Comparison of LSTM-SSP and LSTM-LTP predictions. A well-trained LSTM achieves accurate short-term predictions, but its performance degrades for long-term prediction when true states are no longer provided as conditioning inputs.}
\label{fig:LSTM_SSP_Long_term}
\vspace{-0.35cm}
\end{figure}

We next evaluate the predictive accuracy of the LSTM-SSP and LSTM-LTP schemes introduced in Section~\ref{3.3:DA+ML}. Recall that, in the LSTM-SSP setting, each one-step-ahead prediction is conditioned on the true system states from previous time steps, whereas in the LSTM-LTP setting, predictions are generated recursively using the model’s own past predictions. Figure~\ref{fig:LSTM_SSP_Long_term} illustrates the state prediction performance for the Lorenz–96 system in the $11$th component (i.e., $x_t^{11}$) under the two approaches. The blue curve represents the true state trajectory obtained from direct simulation of the Lorenz-96 system, the black curve shows the LSTM-SSP predictions, and the red curve shows the LSTM-LTP predictions. As shown, the LSTM-SSP predictions remain close to the true trajectory, since this approach conditions each one-step-ahead prediction on the true system states. In contrast, under the LSTM-LTP approach, prediction errors accumulate recursively over time, resulting in a progressive divergence from the true trajectory. In practical applications, however, true system states are typically unavailable, and AI-based forecasts therefore more closely reflect the performance of the LSTM-LTP setting.

\subsubsection{Experimental setup for data assimilation} \label{sec:LSTM_EnSF}

The above results demonstrate that effective uncertainty reduction is essential for improving AI-based predictive performance. In this work, we leverage the generative AI-enabled EnSF to systematically reduce predictive uncertainty through the data assimilation framework.
The data used for uncertainty reduction are obtained through observational operators, and we consider a challenging yet practical setting in which the available observations are only partial and corrupted by noise. Specifically, at each data assimilation step, the true state is observed through a nonlinear operator $h(\cdot)$ (applied component-wise) and is perturbed by Gaussian observational noise. In this work, we adopt the arctangent function, $h(\cdot) = \arctan(\cdot)$, whose bounded output $\left(-\frac{\pi}{2}, \frac{\pi}{2}\right)$ can strongly distort the observations when the magnitude of the true state exceeds this range. 
Consequently, the observation process is as
\begin{equation}\label{eqn:LSTM_Lorenz96_obs_function}
    Y_n = h(X_n) + E^O = \arctan(X_n) + E^O, \quad n = 1, 2, \cdots, 
\end{equation}
where $X_n$ denotes the true state at time step $t_n$ and the observational noise is modeled as Gaussian, e.g., $E^O \sim N(0, 0.01^2I_{20})$.

\begin{figure*}[h!]
\vspace{-0.5cm}
\centering \includegraphics[width=0.7\linewidth]{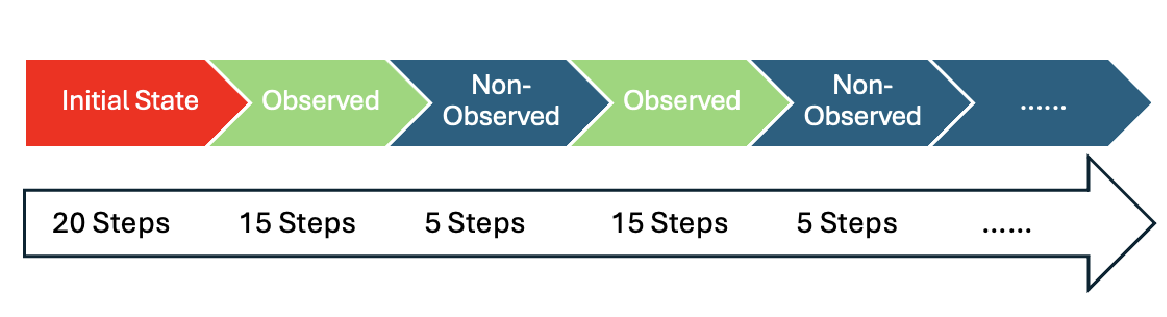}
\caption{Schematic of the data assimilation procedure. After the initial 20 simulation steps, observations are available for 15 out of every 20 time steps.}
\label{fig:DA_Process}
\vspace{-0.5cm}
\end{figure*}

\vspace{0.5em}
The primary goal of our numerical experiments is to demonstrate that EnSF-based data assimilation can effectively reduce predictive uncertainty and thereby improve prediction accuracy when AI models are less accurate or unreliable. To this end, we adopt the LSTM-LTP framework, which enables long-term prediction, as the forward dynamical model, e.g., function $f$ in Eq.~\eqref{Sec3:DA_Framework}, within the data assimilation framework.
To initiate the EnSF, we generate an ensemble of $N=100$  members drawn from the prescribed initial state distribution. Each ensemble member is then propagated forward by integrating the Lorenz–96 system, yielding the state trajectories over the first 19 time steps, $\{X_{1}^{i}, \hdots, X_{19}^{i}\}^{20}_{i=1}$. In this experiment, we adopt a 20-step observation cycle, in which observations are assimilated during the first 15 steps of each cycle, followed by 5 steps without observations. The cycle repeats until the time step $2500$ is reached (see Figure~\ref{fig:DA_Process} for an illustration). 

\subsubsection{Numerical results for a sufficiently trained AI model}
In our first testing scenario, we assume that the LSTM model is \emph{sufficiently trained} using \emph{high-quality, accurate training data} that faithfully represent the underlying dynamical system, which in this case is the Lorenz–96 model.

\begin{figure}[htb]
\centering
\begin{minipage}{0.75\textwidth}
    \centering
    \includegraphics[scale = 0.5]{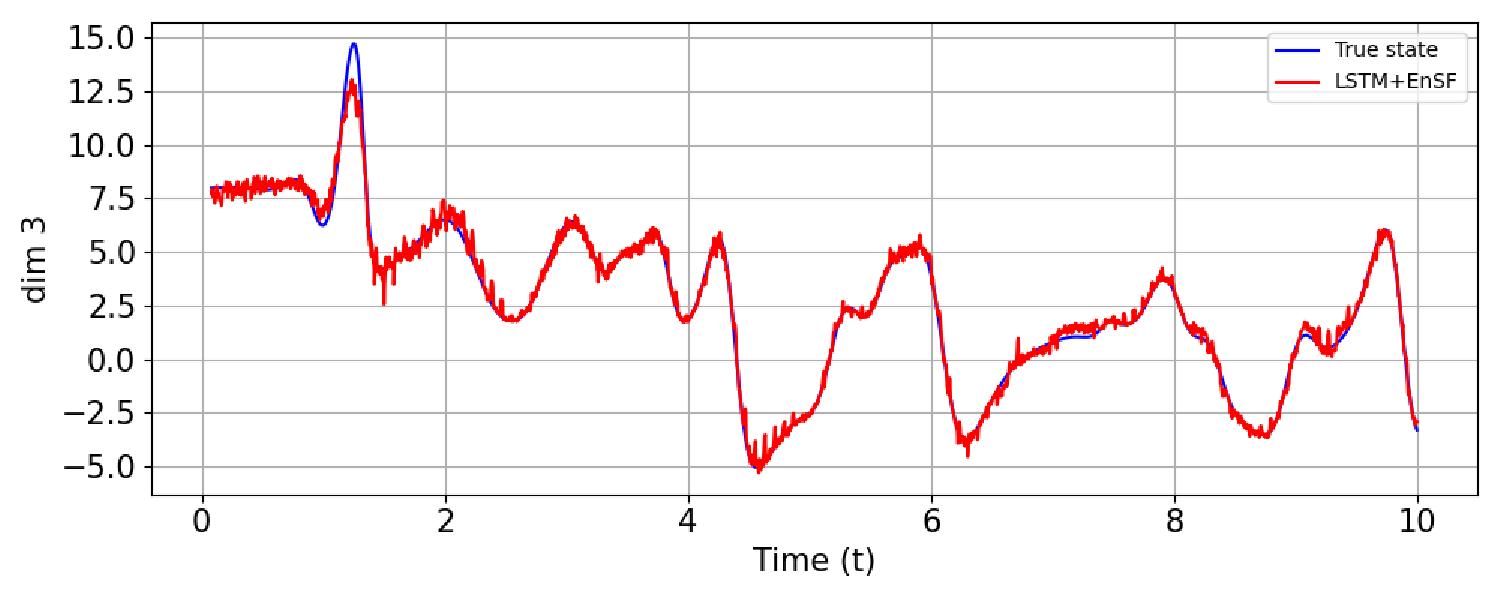}
\end{minipage}
\begin{minipage}{0.75\textwidth}
    \centering
    \includegraphics[scale = 0.5]{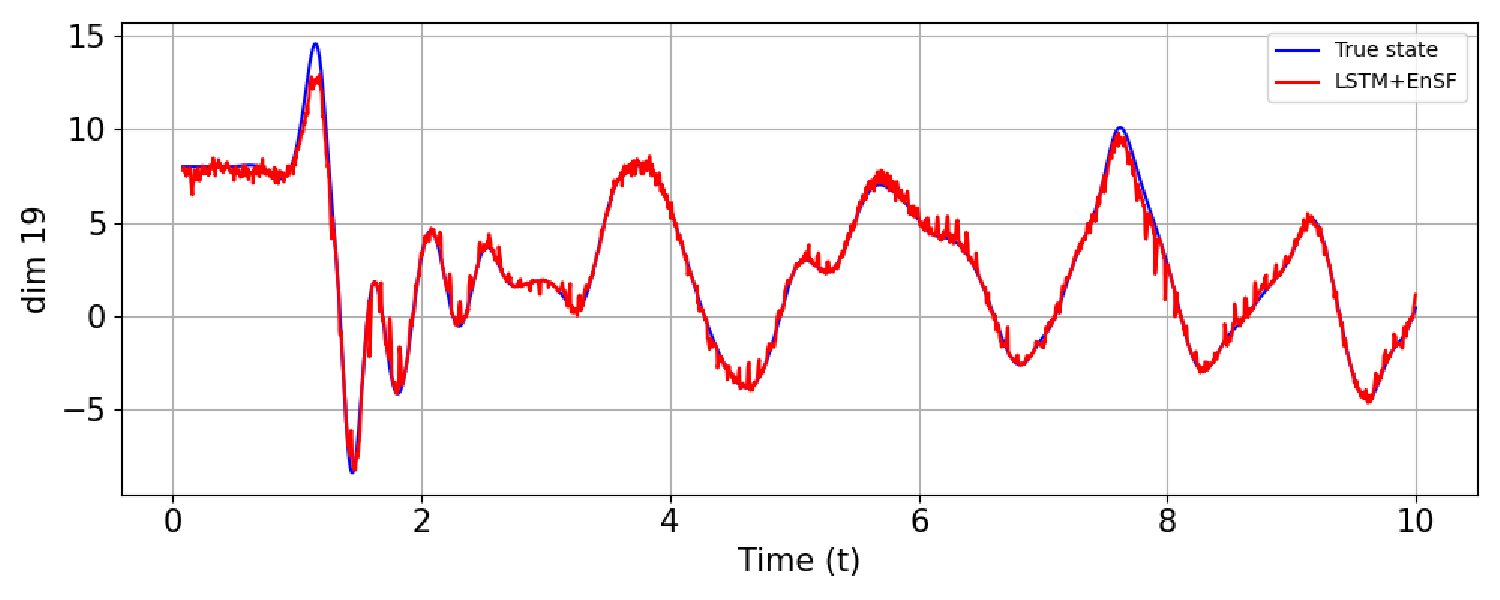}
\end{minipage}
\caption{Performance of the EnSF in calibrating state predictions. (Top) State dimension 3. (Bottom) State dimension 19. The blue curves represent the true states, while the red curves show the LSTM predictions calibrated by the EnSF. The LSTM + EnSF framework achieves accurate and stable state prediction across both dimensions.}
\label{fig:LSTM_SA}
\vspace{-0.45cm}
\end{figure}

\begin{figure}[htb]
\centering
\includegraphics[scale = 0.6]{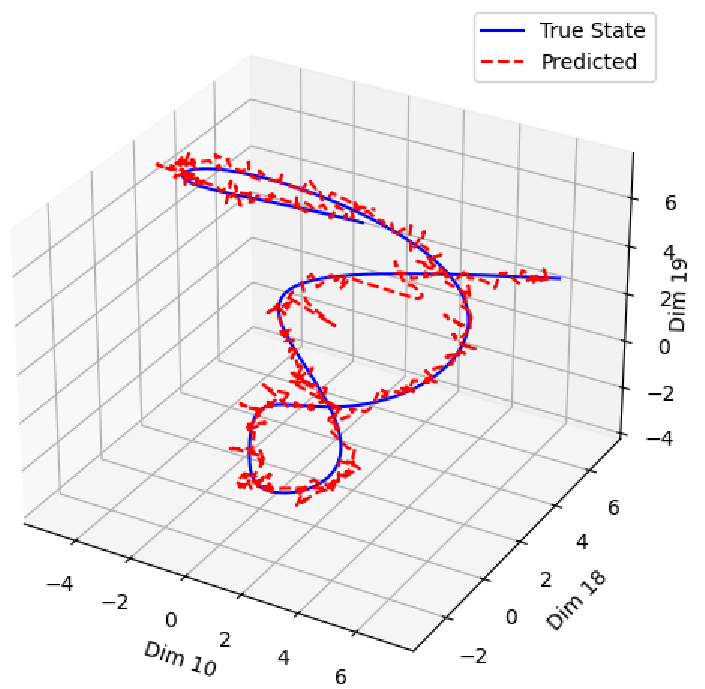}
\caption{3D visualization of state estimation performance over 500 time steps using three randomly selected state components (dimensions 10, 18, and 19).}
\label{fig:LSTM_SA_3D}
\vspace{-0.5cm}
\end{figure}

We begin by considering a fully observed setting, in which all 20 state dimensions are available. The true state trajectory of the Lorenz–96 system in Eq.~\eqref{eqn:Lorenz96_Definition} is obtained via direct numerical simulation of the model. Figure~\ref{fig:LSTM_SA} shows the performance of our method in recovering the true state trajectory for representative state dimensions 3 and 19. To better illustrate the accuracy of the EnSF-calibrated state trajectories, we also present a three-dimensional (3D) visualization of the state estimation performance over 500 time steps in Figure~\ref{fig:LSTM_SA_3D}. In comparison with the pure long-term predictions of the LSTM-LTP model shown in Figure~\ref{fig:LSTM_SSP_Long_term}, the EnSF-enhanced model consistently yields accurate predictions, resulting in a more reliable and accurate LSTM-LTP–based solver.


\begin{figure}[htb]
\centering
\includegraphics[scale = 0.5]{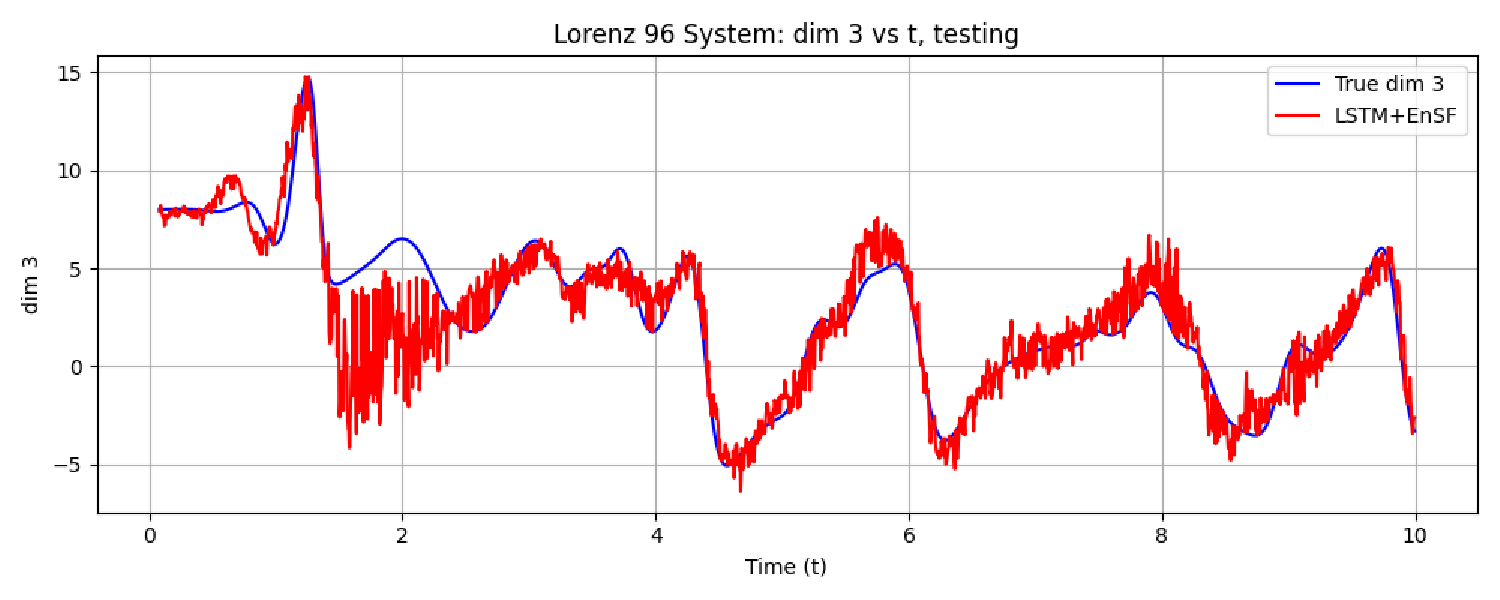}
\caption{LSTM-LTP prediction with EnSF along dimension 3 with partial observational data.}
\label{fig:SA_Partial_obs_12}
\vspace{-0.35cm}
\end{figure}

We next consider a more challenging scenario in which the state is only \emph{sparsely observed}, with only a subset of state components available at each assimilation time. This setting is practically motivated, as it is typically challenging for sensor networks to provide measurements of all target variables at all locations simultaneously, which leads to inherent spatiotemporal sparsity in operational data collection. In this experiment, we still use a $20$-step cycle consisting of $15$ consecutive data assimilation steps followed by $5$ forecast-only steps (without data assimilation). \emph{During each of the $15$ data assimilation steps, we assume that we can only observe $5$ out of the $20$ target state components, and the observable states are chosen randomly.} 
We plot the true and predicted trajectories of a randomly selected component in Figure~\ref{fig:SA_Partial_obs_12}.
Compared with Figure~\ref{fig:LSTM_SA}, the predicted trajectory in this experiment is noisier and exhibits more frequent oscillations, which is as expected due to the sparsity of the observational data. Nevertheless, the EnSF effectively reduces the error accumulated during the LSTM forecast, thereby keeping the predicted state close to the ground truth.

\begin{figure}[htb]
\centering
\includegraphics[scale = 0.4]{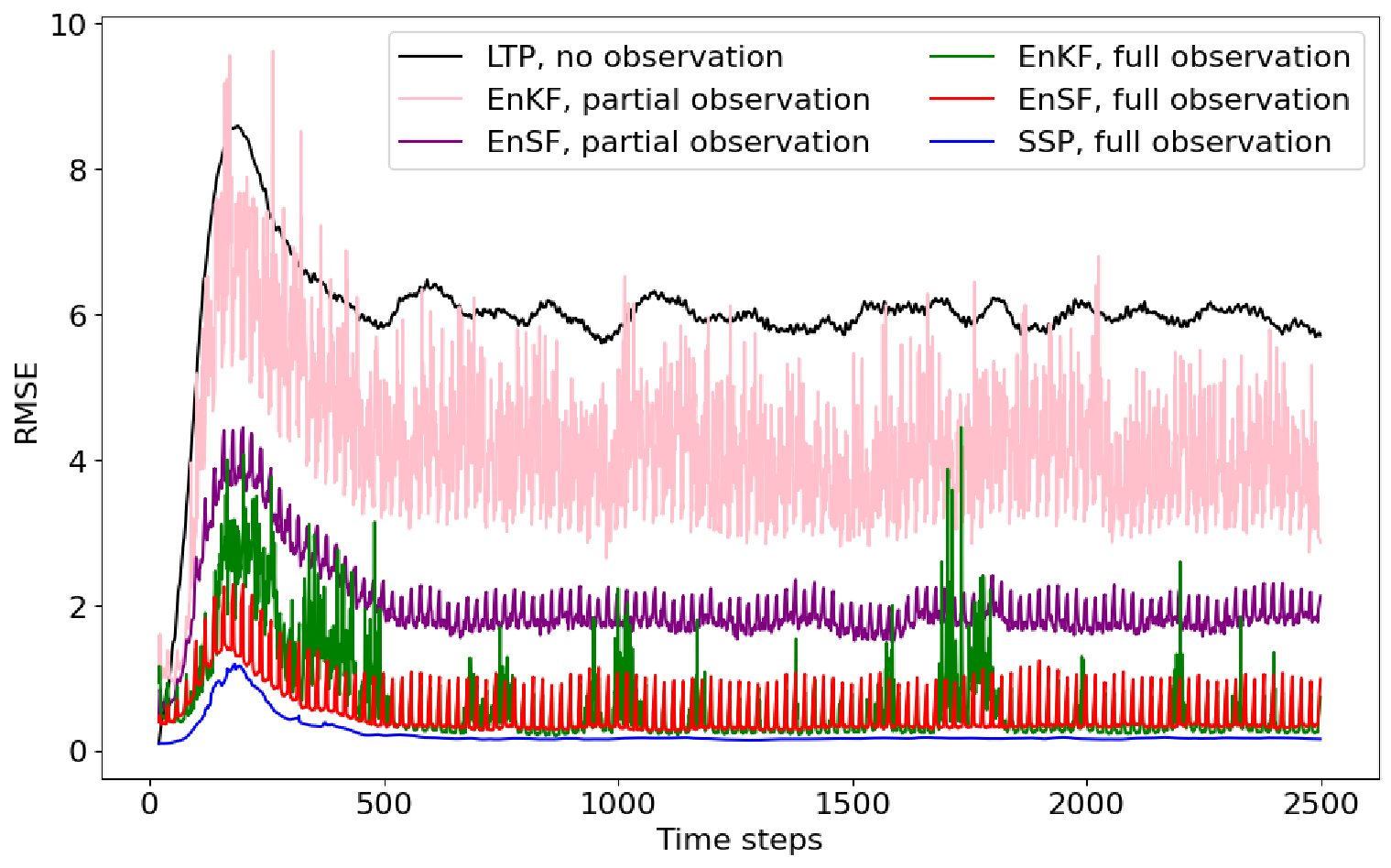}
\caption{RMSE over 40 repeated tests for different prediction settings. The blue curve corresponds to single-step prediction conditioned on \emph{true states}. The red and green curves show EnSF- and EnKF-enhanced state predictions with full observations, respectively. The purple and pink curves show EnSF- and EnKF-enhanced predictions with partial observations, respectively. The black curve represents the long-term baseline prediction of the LSTM model without observational updates.}
\label{fig:LSTM_Lorenz96_RMSE_all}
\vspace{-0.35cm}
\end{figure}

To further assess the effectiveness of EnSF, we repeat the above experiment 
$40$ times using independent state trajectories and report the resulting root mean square errors (RMSE). We also compare EnSF with the ensemble Kalman filter (EnKF), a widely used benchmark method for operational data assimilation. The corresponding RMSE curves are shown in Figure~\ref{fig:LSTM_Lorenz96_RMSE_all}.
As expected, the short-term predictions (SSP, blue) achieve the lowest RMSE, since they are conditioned on the true states, whereas the AI-based long-term predictions (LTP, black) exhibit the highest RMSE due to error accumulation in the absence of observational correction, consistent with the results shown in Figure~\ref{fig:LSTM_SSP_Long_term}. Under full observations, both EnSF-enhanced (red) and EnKF-enhanced (green) predictions closely track the SSP baseline; however, EnKF displays noticeably larger oscillations than EnSF. Under sparse observations, EnSF (purple) consistently outperforms EnKF (pink), indicating superior robustness to limited and nonlinear measurement information. Finally, all RMSE curves display increased errors near the 200-timestep region, which is attributable to the strongly chaotic dynamics of the Lorenz–96 system, a phenomenon commonly observed in long-horizon simulations of this model.

\subsubsection{Numerical results for insufficiently trained LSTM models}

Up to this point, we have considered an LSTM model trained using sufficient high-quality data. In practice, however, available training data may be noisy or limited in quantity, resulting in a poorly trained and therefore unreliable surrogate model (as illustrated, for example, in Figure~\ref{fig:LSTM_training_MSE_total}). In what follows, we consider two training scenarios. The first AI model, denoted by SI (sufficient, but inaccurate data), is trained on noisy data for 1500 epochs, as described in Section~\ref{sec:LSTM_model_training}. The second model, denoted by IA (insufficient, but accurate data), is trained on high-quality, noise-free clean data but for only 100 epochs. We then compare the performance of EnSF- and EnKF-enhanced AI predictions under these two settings.

\begin{figure}[htb]
\centering
\begin{minipage}{0.5\textwidth}
\centering
\includegraphics[scale = 0.325]{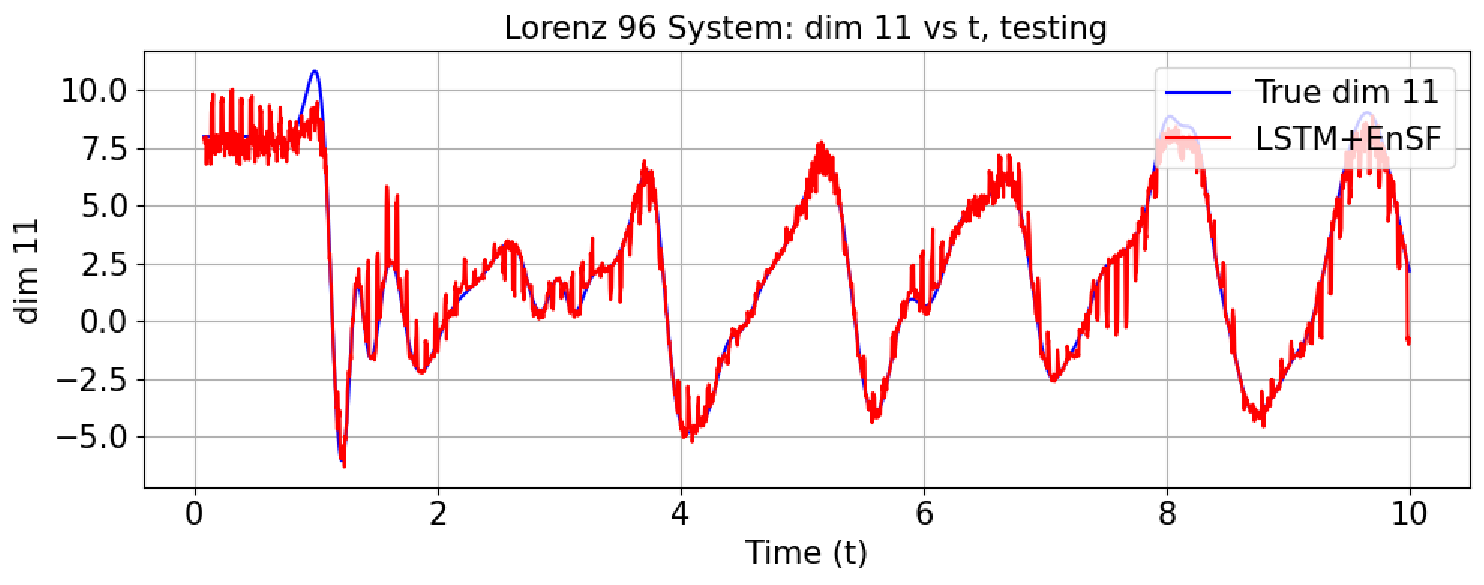}
\label{fig:LSTM_SI_12}
\end{minipage}%
\begin{minipage}{0.5\textwidth}
\centering
\includegraphics[scale = 0.325]{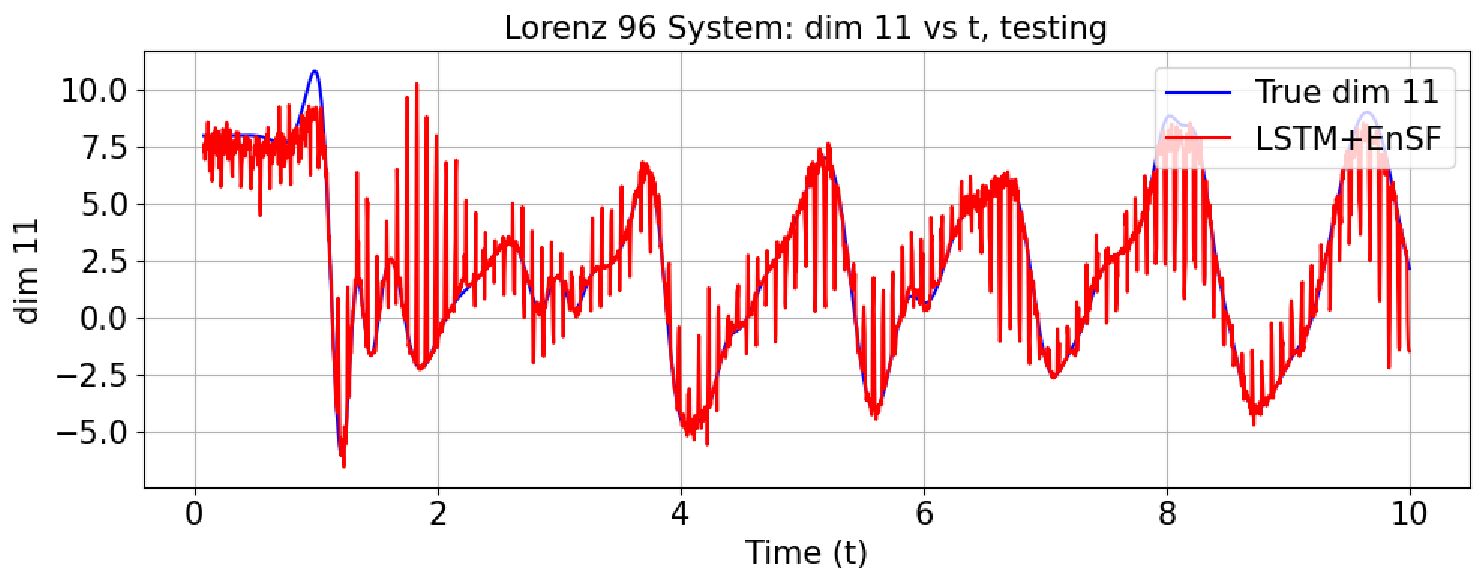}
\label{fig:LSTM_IA_12}
\end{minipage}%

\begin{minipage}{0.5 \textwidth}
\centering
\includegraphics[scale = 0.35]{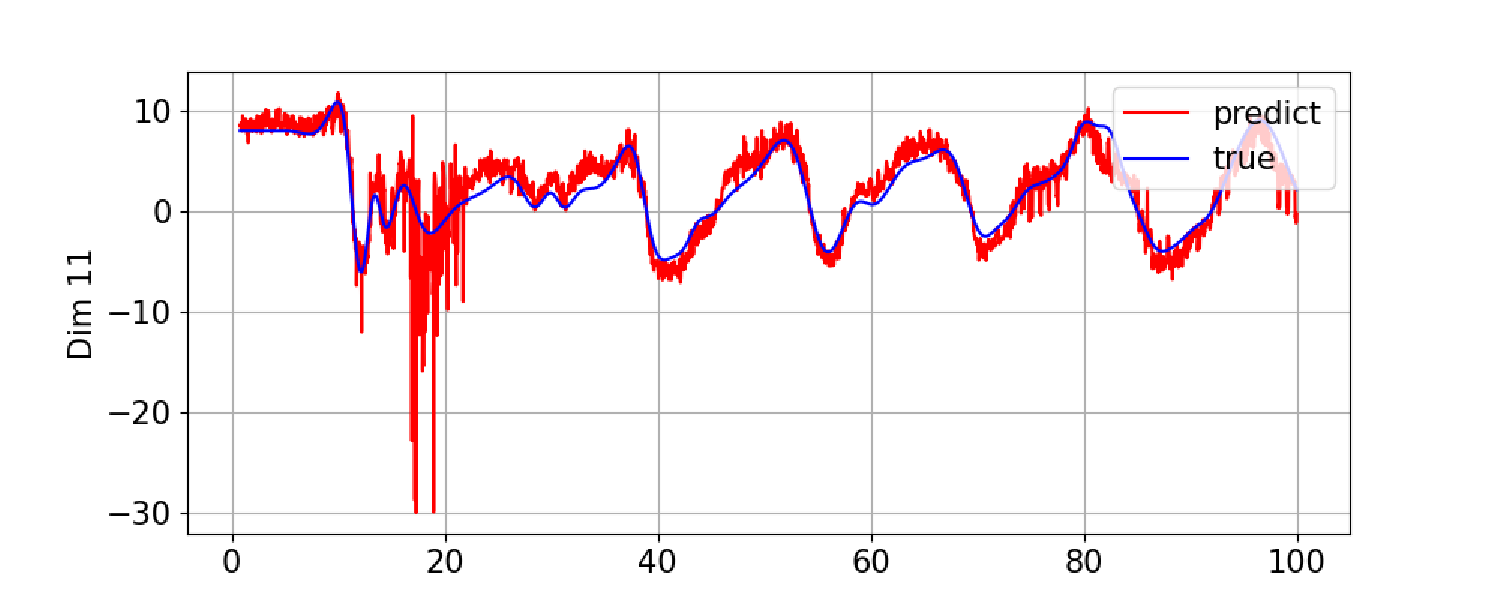}      \label{fig:LSTM_SI_EnKF_12}
\end{minipage}%
\begin{minipage}{0.5 \textwidth}
\centering
\includegraphics[scale = 0.35]{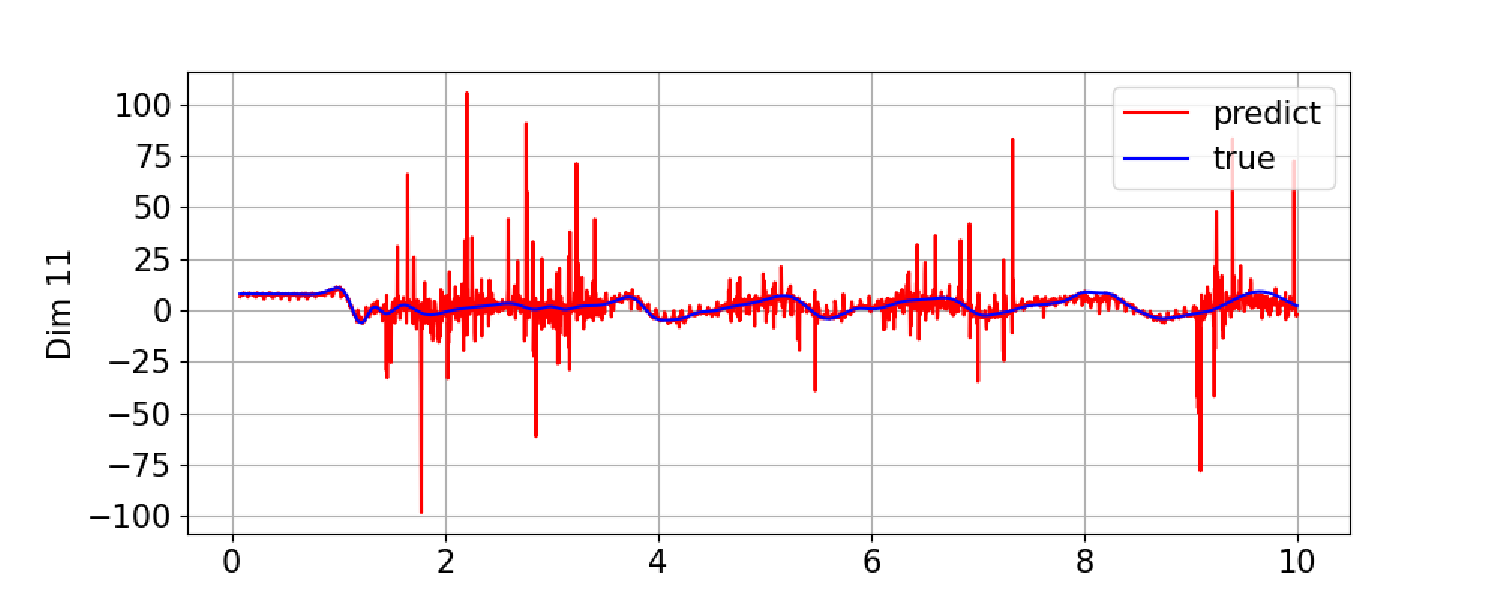}      \label{fig:LSTM_IA_EnKF_12}
\end{minipage}
\vspace{-0.4cm}
\caption{Comparison of the predicted state trajectory along dimension 11 enhanced by EnSF and EnKF under imperfectly trained LSTM models. Top left: SI with EnSF. Top right: IA with EnSF. Bottom left: SI with EnKF. Bottom right: IA with EnKF. Overall, EnKF-enhanced state estimates exhibit larger errors and more pronounced oscillations compared with those enhanced by EnSF. }  \label{fig:LSTM_SI_IA_EnSF_EnKF}
\vspace{-0.35cm}
\end{figure}

The predicted state trajectories enhanced by EnSF and EnKF are shown in Figure~\ref{fig:LSTM_SI_IA_EnSF_EnKF}. Overall, EnSF outperforms EnKF and yields more accurate updated state estimates for both models. In particular, SI-EnSF and IA-EnSF track the true trajectory more closely across the prediction horizon than their EnKF counterparts, although IA-EnSF is more oscillatory due to limited training data. In contrast, SI-EnKF and IA-EnKF exhibit clear instability, producing spuriously large amplitudes and substantial deviations from the true state.

\begin{figure}[htb]
\centering
\includegraphics[scale = 0.4]{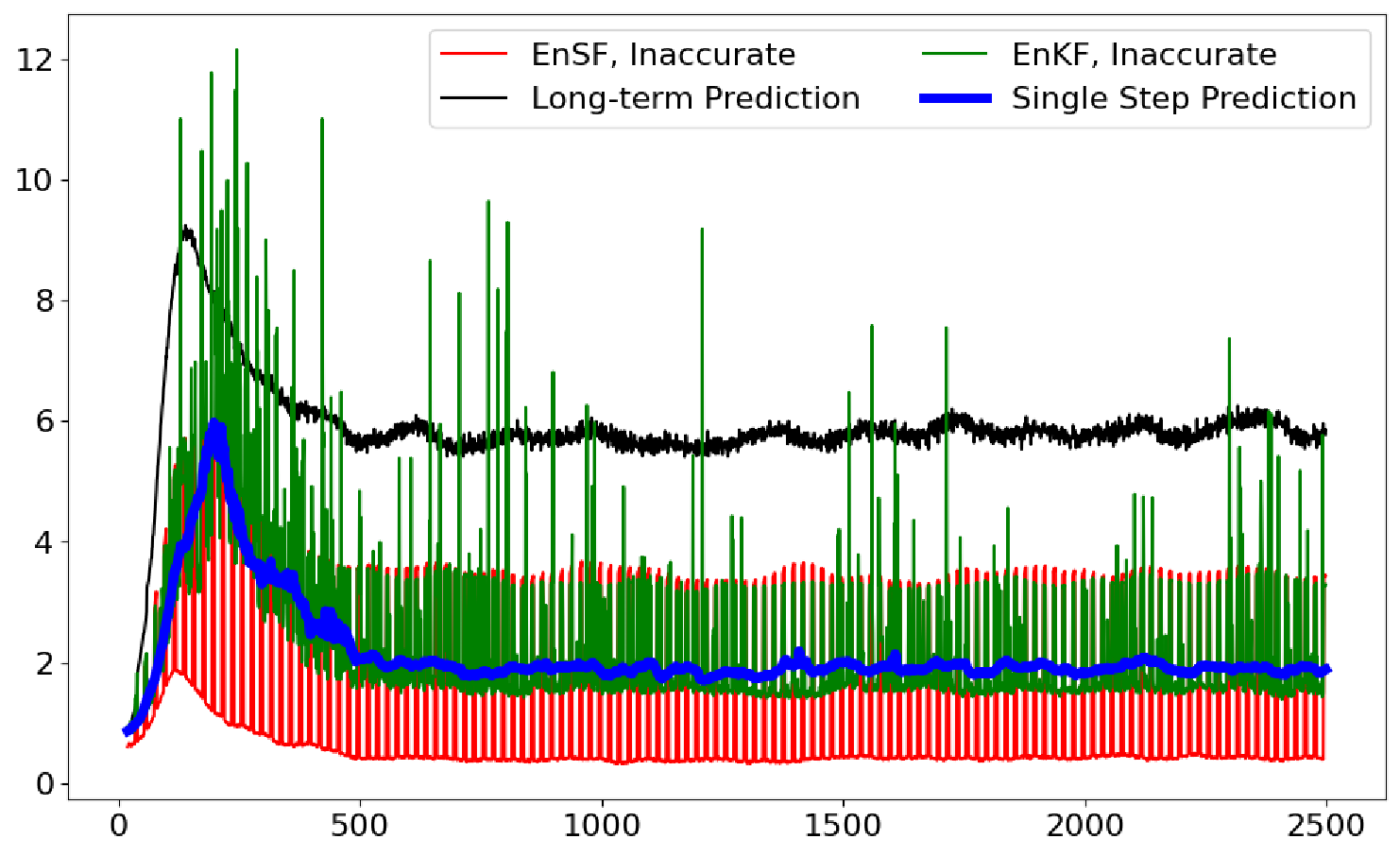}
\caption{RMSEs for SI model with 40 repeated tests.}
\label{fig:LSTM_Lorenz96_RMSE_all_inaccurate}
\end{figure}
\begin{figure}[htb]
\centering
\includegraphics[scale = 0.4]{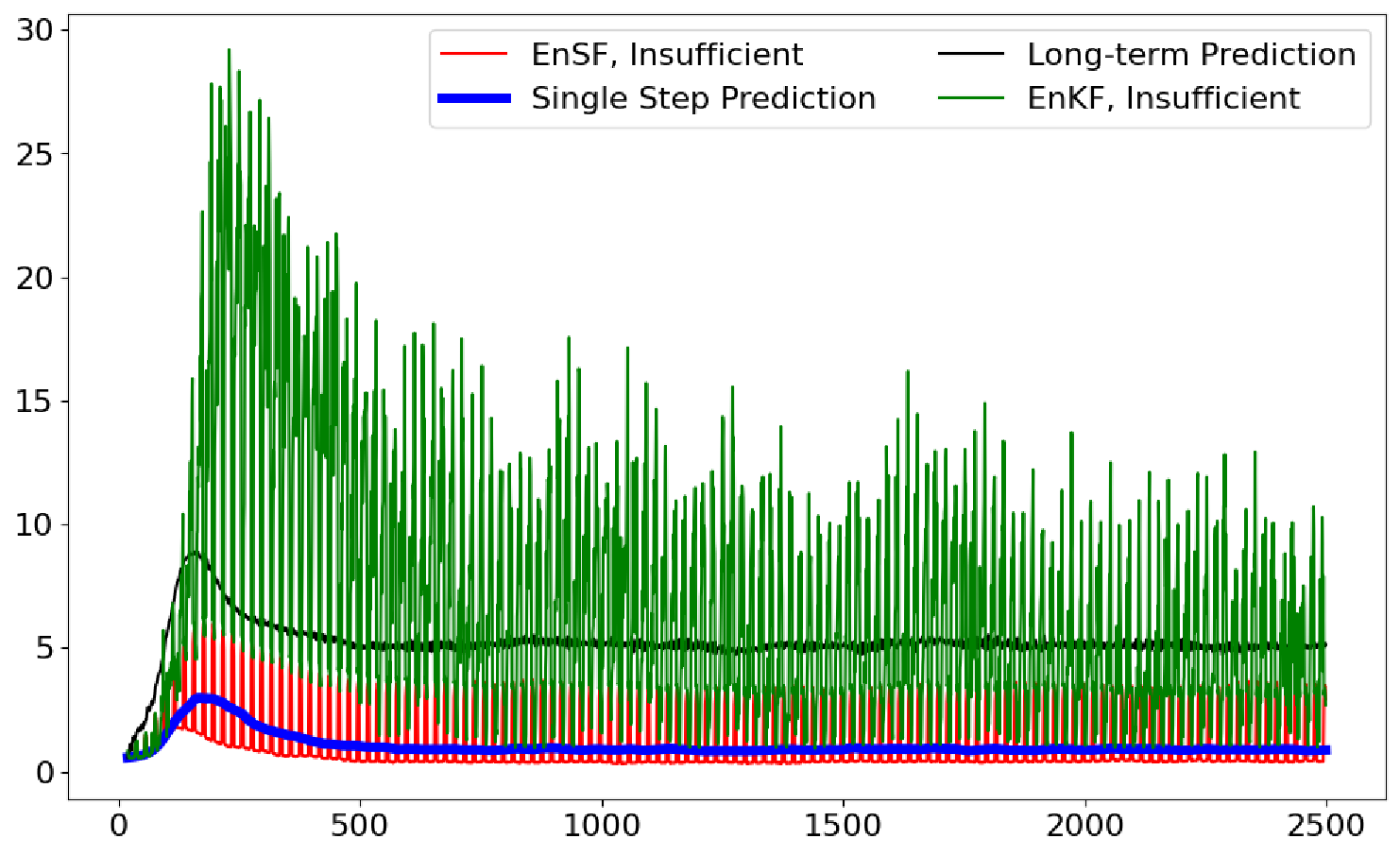}
\caption{RMSEs for IA model with 40 repeated tests.}  \label{fig:LSTM_Lorenz96_RMSE_all_insufficient}
\vspace{-0.45cm}
\end{figure}

We further validate the above results by performing $40$ repeated tests. The corresponding RMSEs are shown in Figures~\ref{fig:LSTM_Lorenz96_RMSE_all_inaccurate} for the SI model and~\ref{fig:LSTM_Lorenz96_RMSE_all_insufficient} for the IA model, together with the single-step prediction (SSP) and long-term prediction (LTP) baselines obtained using the same imperfect forward models.

Across both figures, EnSF consistently yields smaller and more stable RMSEs than EnKF.
Both methods display pronounced RMSE peaks that align with the adopted 20-step observations cycle, which also has been shown in the previous section. During the 5 forecast-only intervals, the state evolution relies entirely on the learned forward model; consequently, model errors can accumulate rapidly and manifest as sharp RMSE spikes. Once data assimilation resumes, both EnSF and EnKF reduce their errors significantly, with EnSF consistently achieving much lower RMSEs across the assimilation steps, indicating a stronger capability to correct imperfect forecasts.

Moreover, although SSP attains the smallest errors in the previous setting with sufficient, high-quality training data (see Figure~\ref{fig:LSTM_Lorenz96_RMSE_all}), this behavior does not persist when the forward model becomes less reliable, as in the SI and IA scenarios.
As shown in Figures~\ref{fig:LSTM_Lorenz96_RMSE_all_inaccurate} and~\ref{fig:LSTM_Lorenz96_RMSE_all_insufficient}, the data assimilation-enhanced state estimates are constructed using observations from the true system. Consequently, they can yield more accurate state reconstructions than the model-based SSP forecasts when the underlying predictive model is inaccurate or insufficiently trained. In these cases, data assimilation effectively corrects the forecast trajectories, leading to improved overall prediction accuracy.

\vspace{1em}
\subsection{Korteweg-De Vries (KdV) equation}\label{4.2:KdV}
To demonstrate the effectiveness of data assimilation in reducing predictive uncertainty in AI models, we next consider using machine learning–based models to solve the following Korteweg–de Vries (KdV) equation:
\vspace{-0.3cm}
\begin{equation}
\label{Eq:KdV}
v_t = -6vv_x -v_{xxx}, \quad 0 \leq x \leq 50, 
\end{equation}
which describes the propagation of shallow-water waves \cite{michalowska2024}. Our initial condition is taken as a superposition of two solitary waves, i.e., $v(x, 0) = \sum\limits^2_{j=1}\frac{1}{2}c_j\cosh^2\left(\frac{\sqrt{c_j}}{2}x-a_j\right)$ 
where $a_j$ controls the initial position of the $j$th wave and $c_j$ controls its amplitude. In this section, we use the R-DeepONet solver described in Section~\ref{2.3:Combination} as the AI surrogate model and study both well-trained and poorly trained cases.

\subsubsection{Training configuration for the R-DeepONet model}
The training set is constructed based on $30$ randomly sampled initial conditions, where each initial condition is parameterized by two independently drawn pairs $\left(c_j, a_j\right)^{2}_{j=1}$. 
The training inputs are constructed using a 10-step history window, consisting of sequences of length 10 sampled at 100 fixed spatial points. 
Unlike the standard DeepONet formulation, our architecture incorporates a recurrent LSTM-based branch network and therefore does not treat time as an explicit input variable. Instead, each training sample is formed by taking the Cartesian product of (i) a length-10 history window of solution snapshots and (ii) the operator output evaluation points. 
Moreover, since the data assimilation propagation satisfies $v_k=f(v_{k-1})$, both $v_k$ and $v_{k-1}$ live in the same spatial state space (i.e., on the same grid). Therefore, we choose the operator output points to coincide with the spatial sensor locations used at future steps.

For the branch network, we employ an LSTM with three hidden layers of width 64, followed by a fully connected layer of width 100. The trunk network consists of three fully connected hidden layers with widths 64, 100, and 100, respectively. This design ensures that the branch and trunk networks produce latent vectors of the same dimension so that their inner product can be used to form the DeepONet output.

We train the operator network sequentially over the set of KdV initial conditions, running 100 epochs for each initial condition before proceeding to the next. We use the Adam optimizer with learning rate $10^{-3}$ and batch size 500. Since operator-learning models can exhibit highly fluctuating training loss due to sensitivity of the architecture \cite{VAnh2024}, learning-rate scheduling is often beneficial; however, because the number of epochs per initial condition is relatively small, we do not apply a learning-rate scheduler in these experiments.

The experimental setup for data assimilation follows the same procedure described in Section~\ref{sec:LSTM_EnSF}. In what follows, we validate the performance of EnSF-enhanced predictions for both sufficiently trained and insufficiently trained models.

\subsubsection{Numerical results for sufficiently trained R-DeepONet model}
We first examine the case that the R-DeepONet model is adequately trained. As discussed in previous sections, under the LTP setting -- where the surrogate model is iteratively rolled out and its own predictions are fed back as inputs -- forecast errors accumulate quickly over time. We therefore apply data assimilation as an uncertainty reduction technique to assess whether the resulting behavior is consistent with that observed in the Lorenz-96 case. Specifically, we consider both EnSF and EnKF and compare their performance against pure LTP without data assimilation.

\begin{figure*}[h!]
\centering
\includegraphics[width=0.95\linewidth]{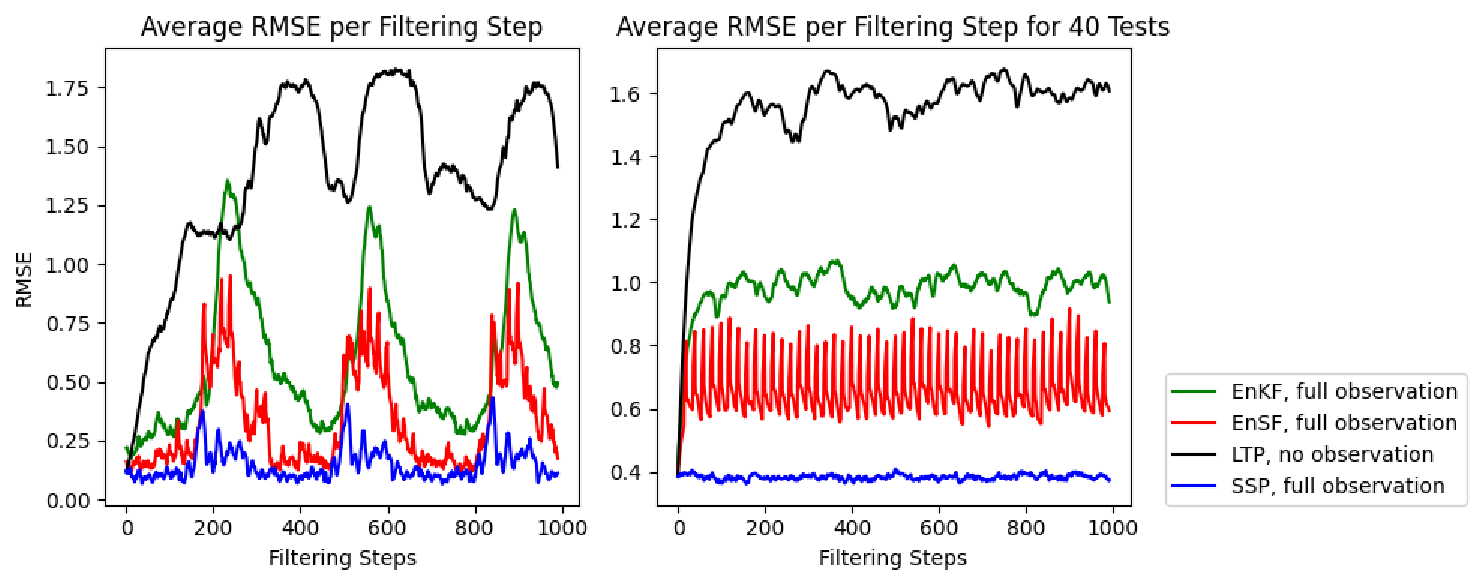}
\caption{RMSE Comparison. (Left) $40$ repeated tests with the same initial condition. (Right) $40$ repeated tests with different initial conditions. The EnSF-enhanced prediction (red) outperforms the EnKF-enhance prediction (green).}
\label{RMSEs_FullTrain_Compare}
\end{figure*}

Figure~\ref{RMSEs_FullTrain_Compare} presents the RMSEs. The left panel shows RMSEs from 40 repeated tests using the same initial condition, while the right panel shows RMSEs from 40 repeated tests with different initial conditions.  Without the data assimilation procedure to reduce uncertainty, the LTP rollout (black) deteriorates rapidly. On the other hand, the data assimilation enhanced LTP predictions improve substantially, with much smaller RMSE that stabilizes quickly. Figure~\ref{RMSEs_FullTrain_Compare} also shows that EnSF-enhanced prediction (red) again outperforms EnKF-enhance prediction (green) and yields lower errors. The SSP baseline (blue), which represents the best achievable prediction obtained by using the true state as model input, attains the lowest RMSEs overall.

\subsubsection{Numerical results for insufficiently trained R-DeepONet model}

We next consider the case where the R-DeepONet model is not well trained. As in the previous section, we study two settings: an \textbf{insufficient-but-accurate-data} (IA) model, trained on noise-free solution trajectories, and a \textbf{sufficient-but-inaccurate-data} (SI) model, trained on solution trajectories that perturbed by Gaussian noise with scale 0.5. 

\begin{figure}[htb]
\centering
\includegraphics[width=0.95\linewidth]{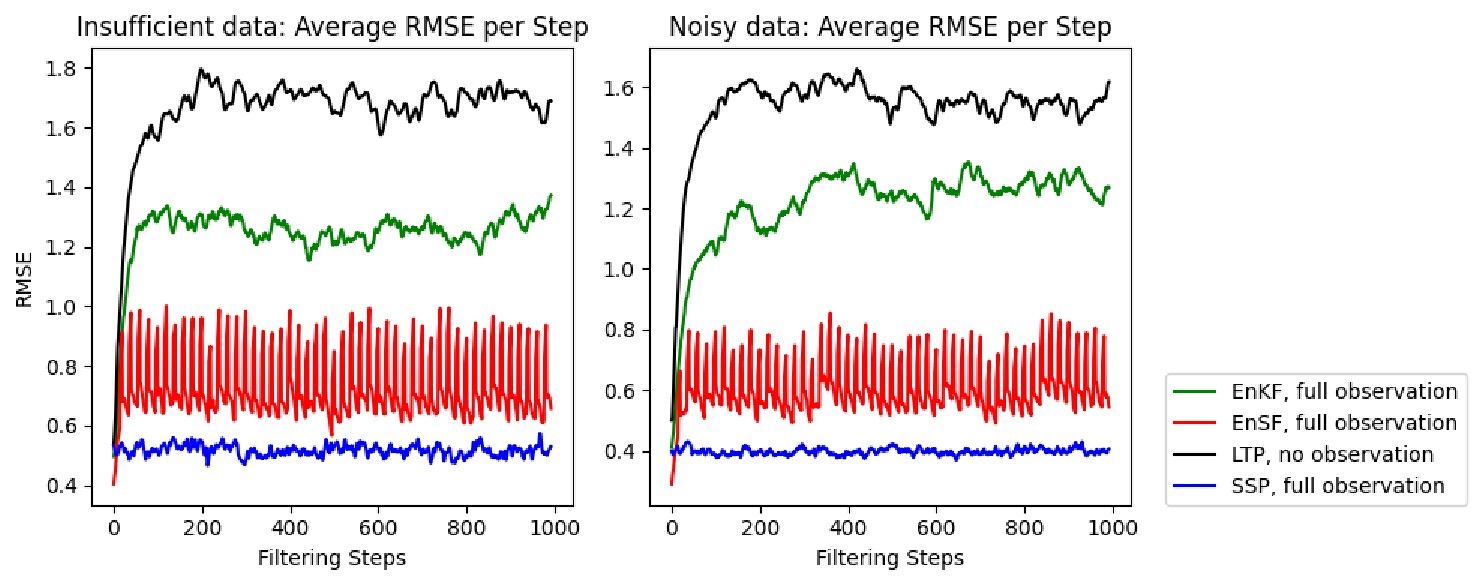}
\caption{Comparison of RMSE for 40 repeated tests with different initial conditions. (Left)) R-DeepONet is trained with insufficient but accurate data. (Right) R-DeepONet is trained with sufficient but noisy data.}
\label{PoorTrain_KdV}
\end{figure}

Figure~\ref{PoorTrain_KdV} presents the RMSE results computed from 40 repeated tests with different initial conditions for the IA model (trained with insufficient but accurate data) and the SI model (trained with sufficient but noise-perturbed data). The results reveal a consistent trend: SSP, as a baseline accuracy demonstration, achieves the lowest error (blue curves), whereas LTP without data assimilation yields the largest RMSE (black curves). When adopting data assimilation as an uncertainty reduction mechanism, the prediction errors of AI models are are substantially reduced. 

Since the EnKF approach relies heavily on the reliability of the underlying model, its performance deteriorates in this insufficiently trained setting where the R-DeepONet model is poorly trained. As shown in the figure, EnSF (red curves) significantly outperforms EnKF (green curves), indicating that EnSF provides a more accurate and robust uncertainty reduction strategy under model inadequacy.

\subsubsection{Numerical experiments with partial observations} 

For the KdV equation example, we also consider the case of partial observations. Specifically, we use the same SA (“sufficient-and-accurate”) model and hyper-parameters as in the previous section, and assume that at each filtering step we observe only a randomly selected half of the spatial grid points.

First, the forward solver produces predictions at all 100 sensor points in the output domain. At each filtering step, however, EnSF receives observations at only 50 randomly selected sensor locations. The data-assimilation update is applied at these 50 observed points, while the recurrent operator network predictions are retained at the remaining 50 unobserved points when constructing the operator input for the next time step.

\begin{figure}[h!]
\centering
\includegraphics[width=0.95\linewidth]{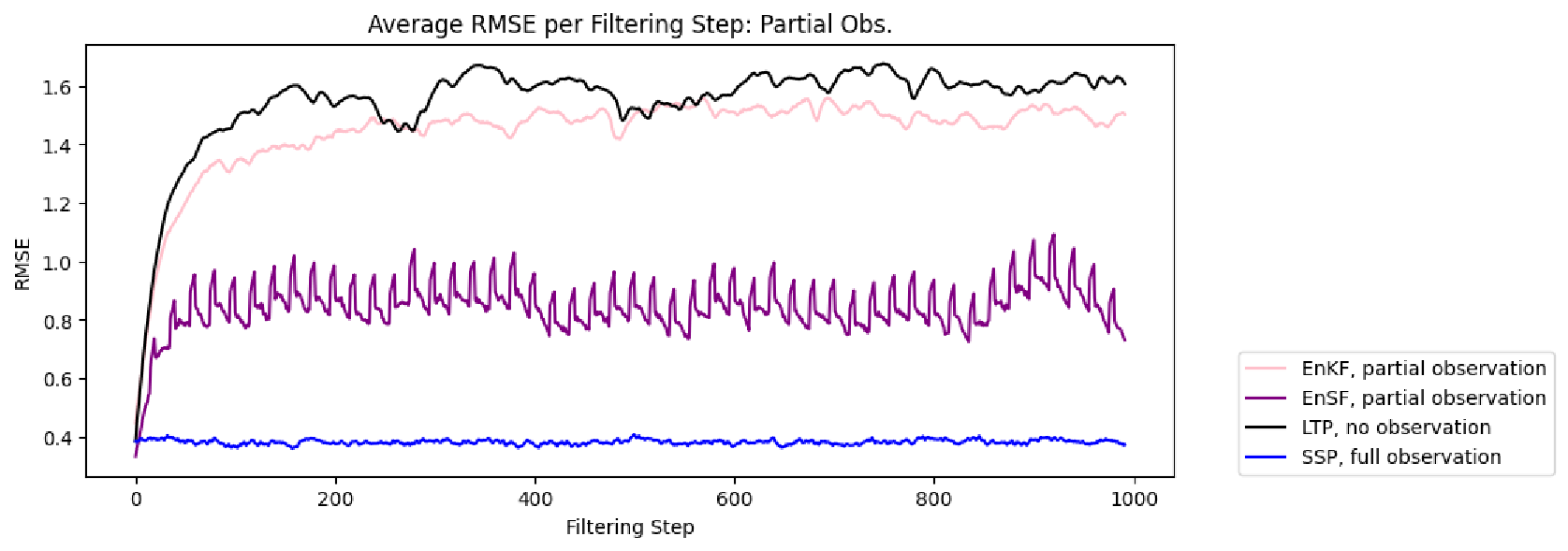}
\caption{Comparison of Average RMSE for 40 KdV Test Surfaces for the R-DeepONet Trained with only $50\%$ of Sensors Locations Observed per Filtering Step}
\label{RMSE_KdVPartial}
\end{figure}

We examine the performance of LTP with EnKF and LTP with EnSF, and report the mean RMSE over 40 test trajectories for the KdV equation from these methods in Figure~\ref{RMSE_KdVPartial}. We also report the corresponding curves from the LTP without data assimilation and the SSP. In this partial-observation setting, EnKF performs poorly, yielding an RMSE comparable to LTP without data assimilation. In contrast, EnSF consistently outperforms EnKF and substantially improves the performance of the LTP. 

\begin{figure}[h!]
\centering
\includegraphics[width=0.95\linewidth]{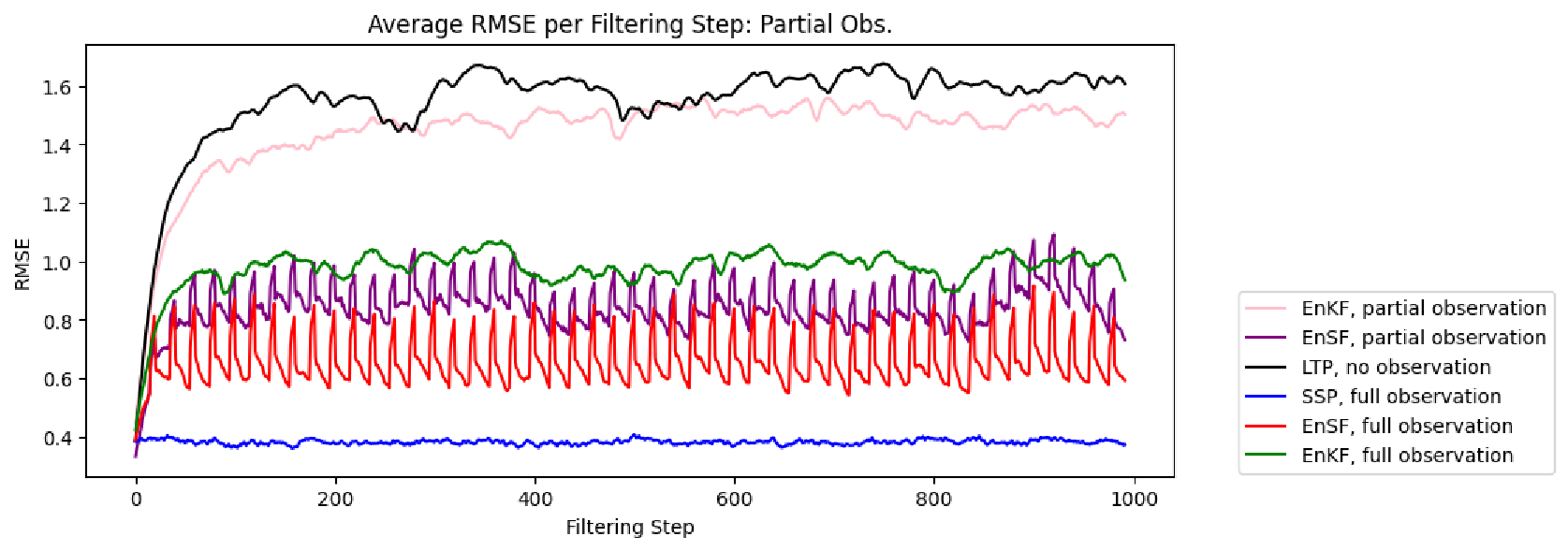}
\caption{Comparison of Average RMSE for 40 KdV Test Surfaces for the R-DeepONet Trained with only $50\%$ of Sensors Locations Observed per Filtering Step}
\label{RMSE_KdVPartial_full}
\end{figure}


\begin{table}[h!]
\centering
\begin{tabular}{|c|c|c|c|c|}
\hline
 &True Data &   Insuff. Data &Noisy Data &Partial\\
 \hline
 \textbf{No Filtering} &1.550 &1.658 & 1.568 & - \\
 \hline
 \hline
\textbf{NN Testing} &0.382 & 0.517& 0.397& -\\
\hline
\hline
\textbf{EnSF}  & 0.517 & 0.731& 0.604& 0.845\\
\hline
EnSF Obs. &0.459  & 0.607& 0.536&0.536 \\
\hline
EnSF Non-Obs. & 0.750 &1.057& 0.886&0.852\\
\hline
\hline
\textbf{EnKF} & 0.977 & 1.246 & 1.215& 1.446\\
\hline
\end{tabular}
\caption{KdV Equation Mean MSE Comparison. EnSF Obs tracks mean MSE only for steps where filtering was applied, while EnSF Non-Obs tracks mean MSE for steps where no filtering was applied.}
\label{tab:OpLearnEnSF}
\vspace{-0.35cm}
\end{table}

Finally, we provide a comprehensive comparison of all methods considered in this section and summarize the mean RMSE in Table~\ref{tab:OpLearnEnSF}. From the first two rows of Table~\ref{tab:OpLearnEnSF}, we can see that SSP yields the lowest error, while LTP without data assimilation performs the worst. Both EnKF and EnSF improve LTP performance, with EnSF remaining more robust than EnKF across all numerical settings. We also report the mean RMSE separately over the 15 data-assimilation steps and the 5 non-assimilation steps within each 20-step observation cycle (denoted by EnSF Obs. and EnSF Non-Obs., respectively). As expected, during the 5 consecutive non-DA steps the RMSE increases noticeably, but these errors are quickly reduced once EnSF updates resume during the online filtering steps. These results highlight the benefit of data assimilation for reducing uncertainty in machine-learning forecasts and demonstrate the effectiveness of the EnSF.

\section{Conclusion}
\label{Sec5:Conclusion}
In this work, we present a hybrid data assimilation-machine learning framework to reduce uncertainty in long-time prediction with machine-learning solvers. We construct forecast models based on the LSTM architecture, which provides a dynamical forward model and enables online filtering with data assimilation. For the assimilation step, we employ the Ensemble Score Filter (EnSF), which is a score-based diffusion model that encodes the underlying physics into the score function using
simulation data and incorporates observational information via a likelihood-based correction in a
reverse-time SDE. This enables iterative refinement of the solution as new data become available, improving predictive accuracy and reducing model uncertainty.

To validate the effectiveness and robustness of our approach, we conduct comprehensive numerical experiments across several scenarios where the learned solver is imperfect. We also compare EnSF with the standard Ensemble Kalman Filter (EnKF). The results demonstrate that EnSF can recover hidden physical states and provide accurate updates even with sparse and noisy observations.

Future directions include extending the proposed framework to stochastic operator learning for SPDEs, and to federated settings where data are distributed across multiple clients and model updates are aggregated to train a global model without sharing raw data.

\bibliographystyle{plain}
\bibliography{References}

@unpublished{Lu2019,
    author = {L. Lu and P. Jin and G. E. Karniadakis},
    title = {Deeponet: Learning nonlinear operators for identifying differential equations based on the universal approximation theorem of operators},
    note = {arXiv preprint, arXiv:1910.03193},
    year = {2019}
}

@article{S-DeepONet,
title = {Sequential Deep Operator Networks (S-DeepONet) for predicting full-field solutions under time-dependent loads},
journal = {Engineering Applications of Artificial Intelligence},
volume = {127},
pages = {107258},
year = {2024},
issn = {0952-1976},
doi = {https://doi.org/10.1016/j.engappai.2023.107258},
url = {https://www.sciencedirect.com/science/article/pii/S0952197623014422},
author = {Junyan He and Shashank Kushwaha and Jaewan Park and Seid Koric and Diab Abueidda and Iwona Jasiuk},
}

@article{Lu2021,
    author = {L. Lu and and G. Pang and P. Jin and Z. Zhang and G. E. Karniadakis},
    title = {Learning nonlinear operators via DeepONet based on the universal approximation theorem of operators.},
    journal = {Nat. Mach. Intell.},
    volume = {3},
    pages = {218-229},
    year = {2021}
}

@unpublished{Kovachki2024,
    author = {N. B. Kovachki and S. Lanthaler and A. M. Stuart},
    title = {Operator learning: Algorithms and analysis},
    note = {arXiv preprint, arXiv:2402.15715},
    year = {2024}
}

@article{cohn1997,
  title={An introduction to estimation theory (gtspecial issueltdata assimilation in meteology and oceanography: theory and practice)},
  author={Cohn, Stephen E},
  journal={Journal of the Meteorological Society of Japan. Ser. II},
  volume={75},
  number={1B},
  pages={257--288},
  year={1997},
  publisher={Meteorological Society of Japan}
}

@article{bao2024ensemble,
title={An ensemble score filter for tracking high-dimensional nonlinear dynamical systems},
author={F. Bao and Z. Zhang and G. Zhang},
journal={Computer Methods in Applied Mechanics and Engineering},
volume={432},
pages={117447},
year={2024},
publisher={Elsevier}
}

@article{BaoF20182,
author = {F. Bao and Y. Cao and W. Zhao},
title = {A backward doubly stochastic differential equation approach for nonlinear filtering problems},
journal = {Comm.  Comp. Phy.},
volume = {23},
year = {201},
pages = {1573-1601}
}

@article{BaoC20142,
author = {F. Bao and Y. Cao and X. Han},
title= {Forward backward doubly stochastic differential equations and optimal Filtering of diffusion processes},
journal = {Communications in Mathematical Sciences},
year = {2020},
volume = {18},
number = {3},
pages = {635-661},
}

@article{Bao_AA20,
	Author = {Bao, F. And Cao, Y. And Yong, J. },
	Journal = {Analysis and Applications},
	Title = {Data informed solution estimation for forward backward stochastic
differential equations},
	Year = {2021},
	volume = {19},
	issue = {03},
	pages = {439-464}
	}

@ARTICLE{EnSF_Scalable_2024,
    AUTHOR = {Yin, J. and Liu, S. and Liang, S. and Bao, F. and Chipilski, H. G. and Lu, D. and Zhang, G.},
    YEAR = {2024},
    TITLE = {{A Scalable Real-Time Data Assimilation Framework for Predicting Turbulent Atmosphere Dynamics}},
    JOURNAL = {The International Conference for High Performance Computing, Networking, Storage, and Analysis (SC24 Workshop)},
    DOI = {10.1109/SCW63240.2024.00009}
 }

@ARTICLE{UF_2026,
    AUTHOR = {Bao, F and Zhang, Z and Zhang, G},
    YEAR = {2026},
    TITLE = {United Filter for Jointly Estimating State and Parameters of Stochastic Dynamical Systems},
    JOURNAL = {arXiv},
    VOLUME = {39},
    Number = {3},
    PAGES = {747-774}
 }

@ARTICLE{IEnSF_2025,
    AUTHOR = {Zhang, Z and Bao, F and Zhang, G},
    YEAR = {2025},
    TITLE = {IEnSF: Iterative Ensemble Score Filter for Reducing Error in Posterior Score Estimation in Nonlinear Data Assimilation},
    JOURNAL = {arXiv},
    VOLUME = {},
    PAGES = {arXiv:2510.20159 }
 }

@inproceedings{michalowska2024,
  title={Neural operator learning for long-time integration in dynamical systems with recurrent neural networks},
  author={Micha{\l}owska, Katarzyna and Goswami, Somdatta and Karniadakis, George Em and Riemer-S{\o}rensen, Signe},
  booktitle={2024 International Joint Conference on Neural Networks (IJCNN)},
  pages={1--8},
  year={2024},
  organization={IEEE}
}

@conference{Kingma2014,
    author = {D. P. Kingma and J. Ba},
    booktitle = { Proceedings of the 3rd International Conference on Learning Representations (ICLR)},
    title = {Adam: A Method for Stochastic Optimization},
    address = {Banff},
    month ={14-16 April},
    year = {2014}
}

@unpublished{VAnh2024,
    author = {V.-A. Le and M. Dik},
    title = {A mathematical analysis of neural operator behaviors},
    note = {arXiv preprint arXiv:2410.21481},
    year = {2024}
}

@article{Lu2022,
    author = {L. Lu and X. Meng and S. Cai and Z. Mao and S. Goswami and Z. Zhang and G. E. Karniadakis},
    title = {A comprehensive and fair comparison of two neural operators (with practical extensions) based on fair data},
    journal = {Computer Methods in Applied Mechanics and Engineering},
    volume = {393},
    year = {2022},
    pages = {114778}
}

@article{Chen1995,
    author = {T. Chen and H. Chen},
    title = {Universal approximation to nonlinear operators by neural networks with arbitrary activation functions and its application to dynamical systems.},
    journal = {IEEE transactions on neural netwrosk},
    volume = {6},
    number = {4},
    year = {1995},
    pages = {911-917}
}

@book{Goodfellow2016,
    title={Deep Learning},
    author={I. Goodfellow and Y. Bengio and A. Courville},
    publisher={MIT Press},
    note={deeplearningbook.org},
    year={2016}
}

@article{bao2024score,
  title={A score-based filter for nonlinear data assimilation},
  author={F. Bao and Z. Zhang and G. Zhang},
  journal={Journal of Computational Physics},
  volume={514},
  pages={113207},
  year={2024},
  publisher={Elsevier}
}

@article{nguyen2004multiple,
  title={Multiple neural networks for a long term time series forecast},
  author={H. H. Nguyen and C. W. Chan},
  journal={Neural Computing \& Applications},
  volume={13},
  number={1},
  pages={90--98},
  year={2004},
  publisher={Springer}
}

@article{raissi2019physics,
  title={Physics-informed neural networks: A deep learning framework for solving forward and inverse problems involving nonlinear partial differential equations},
  author={M. Raissi and P. Perdikaris and G. E. Karniadakis},
  journal={Journal of Computational physics},
  volume={378},
  pages={686--707},
  year={2019},
  publisher={Elsevier}
}

@inproceedings{mao2024optimizing,
  title={Optimizing Modified Barium Swallow Exam Workflow: Automating Pre-Analysis Video Sorting in Swallowing Function Assessment},
  author={S. Mao and A. M. Naser and S. N. Buoy and K. K. Brock and K. A. Hutcheson},
  booktitle={2024 46th Annual International Conference of the IEEE Engineering in Medicine and Biology Society (EMBC)},
  pages={1--4},
  year={2024},
  organization={IEEE}
}

@article{hsieh1998applying,
  title={Applying neural network models to prediction and data analysis in meteorology and oceanography},
  author={W. W. Hsieh and B. Tang},
  journal={Bulletin of the American Meteorological Society},
  volume={79},
  number={9},
  pages={1855--1870},
  year={1998},
  publisher={American Meteorological Society}
}

@article{tsubaki2019compound,
  title={Compound--protein interaction prediction with end-to-end learning of neural networks for graphs and sequences},
  author={M. Tsubaki and K. Tomii and J. Sese},
  journal={Bioinformatics},
  volume={35},
  number={2},
  pages={309--318},
  year={2019},
  publisher={Oxford University Press}
}

@article{greff2016lstm,
  title={LSTM: A search space odyssey},
  author={K. Greff and R. K. Srivastava and J. Koutn{\'\i}k and B. R. Steunebrink and J. Schmidhuber},
  journal={IEEE transactions on neural networks and learning systems},
  volume={28},
  number={10},
  pages={2222--2232},
  year={2016},
  publisher={IEEE}
}

@article{Pathak2018,
  title = {Model-Free Prediction of Large Spatiotemporally Chaotic Systems from Data: A Reservoir Computing Approach},
  author = {J. Pathak and B. Hunt and M. Girvan and Z. Lu and E. Ott},
  journal = {Phys. Rev. Lett.},
  volume = {120},
  issue = {2},
  pages = {024102},
  numpages = {5},
  year = {2018},
  month = {1},
  publisher = {American Physical Society},
  url = {https://link.aps.org/doi/10.1103/PhysRevLett.120.024102}
}

@Article{Allen2025,
author={A. Allen, Anna
and S. Markou
and W. Tebbutt
and J. Requeima
and W. P. Bruinsma
and T. R. Andersson
and M. Herzog
and N. D. Lane
and M .Chantry
and J. S. Hosking
and R. E. Turner},
title={End-to-end data-driven weather prediction},
journal={Nature},
year={2025},
month={5},
day={01},
volume={641},
number={8065},
pages={1172-1179},
issn={1476-4687},
url={https://doi.org/10.1038/s41586-025-08897-0}
}

@Article{Wang2024,
author={M. Wang
and J. Li},
title={Interpretable predictions of chaotic dynamical systems using dynamical system deep learning},
journal={Scientific Reports},
year={2024},
month={2},
day={07},
volume={14},
number={1},
pages={3143},
issn={2045-2322},
url={https://doi.org/10.1038/s41598-024-53169-y}
}

@inproceedings{lorenz1996predictability,
  title={Predictability: A problem partly solved},
  author={E. N. Lorenz},
  booktitle={Proc. Seminar on predictability},
  volume={1},
  pages={1--18},
  year={1996},
  organization={Reading}
}

@article{van2019particle,
  title={Particle filters for high-dimensional geoscience applications: A review},
  author={P. J. Van Leeuwen and H. R. K{\"u}nsch and L. Nerger and R. Potthast and S. Reich},
  journal={Quarterly Journal of the Royal Meteorological Society},
  volume={145},
  number={723},
  pages={2335--2365},
  year={2019},
  publisher={Wiley Online Library}
}

@inproceedings{sohl2015deep,
  title={Deep unsupervised learning using nonequilibrium thermodynamics},
  author={Sohl-Dickstein, Jascha and Weiss, Eric and Maheswaranathan, Niru and Ganguli, Surya},
  booktitle={International conference on machine learning},
  pages={2256--2265},
  year={2015},
  organization={pmlr}
}

@article{snyder2015performance,
  title={Performance bounds for particle filters using the optimal proposal},
  author={Snyder, Chris and Bengtsson, Thomas and Morzfeld, Mathias},
  journal={Monthly Weather Review},
  volume={143},
  number={11},
  pages={4750--4761},
  year={2015}
}

@article{law2015,
  title={Data assimilation},
  author={Law, Kody and Stuart, Andrew and Zygalakis, Kostas},
  journal={Cham, Switzerland: Springer},
  volume={214},
  pages={52},
  year={2015},
  publisher={Springer}
}

@article{bao2025nonlinear,
title={Nonlinear Ensemble Filtering with Diffusion Models: Application to the Surface Quasigeostrophic Dynamics},
author={F. Bao and H. G. Chipilski and S. Liang and G. Zhang and J. S. Whitaker},
journal={Monthly Weather Review},
volume={153},
number={7},
pages={1155-1169},
year={2025},
publisher={American Meteorological Society}
}

@article{Bao2025a,
title = {Ensemble score filter with image inpainting for data assimilation in tracking surface quasi-geostrophic dynamics with partial observations},
journal = {	arXiv:2501.12419 },
author = {S. Liang and H. Tran and F. Bao and H. G. Chipilski and P. Jan van Leeuwen and G. Zhang}
}

@article{archibald2022kernel,
  title={A kernel learning method for backward sde filter},
  author={R. Archibald and F. Bao},
  journal={arXiv preprint arXiv:2201.10600},
  year={2022}
}

@article{hochreiter1997long,
  title={Long short-term memory},
  author={Hochreiter, Sepp and Schmidhuber, J{\"u}rgen},
  journal={Neural computation},
  volume={9},
  number={8},
  pages={1735--1780},
  year={1997},
  publisher={MIT press}
}

@article{gers2000learning,
  title={Learning to forget: Continual prediction with LSTM},
  author={Gers, Felix A and Schmidhuber, J{\"u}rgen and Cummins, Fred},
  journal={Neural computation},
  volume={12},
  number={10},
  pages={2451--2471},
  year={2000},
  publisher={MIT press}
}

@article{sherstinsky2020fundamentals,
  title={Fundamentals of recurrent neural network (RNN) and long short-term memory (LSTM) network},
  author={Sherstinsky, Alex},
  journal={Physica D: Nonlinear Phenomena},
  volume={404},
  pages={132306},
  year={2020},
  publisher={Elsevier}
}

@article{yu2019review,
  title={A review of recurrent neural networks: LSTM cells and network architectures},
  author={Yu, Yong and Si, Xiaosheng and Hu, Changhua and Zhang, Jianxun},
  journal={Neural computation},
  volume={31},
  number={7},
  pages={1235--1270},
  year={2019},
  publisher={MIT Press One Rogers Street, Cambridge, MA 02142-1209, USA journals-info~…}
}

@book{burkhart2019discriminative,
  title={A discriminative approach to Bayesian filtering with applications to human neural decoding},
  author={Burkhart, Michael C},
  year={2019},
  publisher={ProQuest Dissertations Publishing}
}

@article{Chapman1928,
    author ={Chapman, Sydney} ,
    title = {On the Brownian displacements and thermal diffusion of grains suspended in a non-uniform fluid.},
    journal = {Proceedings of the Royal Society of London. Series A, Containing Papers of a Mathematical and Physical Character},
    year = {1928},
    pages ={34--54}
}

@article{Kolmogoroff1931,
    author ={Kolmogoroff, Andrei} ,
    title = {Über die analytischen Methoden in der Wahrscheinlichkeitsrechnung},
    journal = {Mathematische Annalen},
    year = {1931},
    pages = {415 -- 458}
}

@article{lacey1998tutorial,
  title={Tutorial: The kalman filter},
  author={Lacey, Tony},
  journal={Georgia Institute of Technology},
  volume={139},
  year={1998}
}

@article{kalman1960,
title={A new approach to linear filtering and prediction problems},
author={R. E. Kalman},
journal = {J. Basic Eng.},
volume = {82},
number = {1},
pages = {35-45},
year={1960}
}

@article{evensen2003ensemble,
  title={The ensemble Kalman filter: Theoretical formulation and practical implementation},
  author={G. Evensen},
  journal={Ocean dynamics},
  volume={53},
  number={4},
  pages={343-367},
  year={2003},
  publisher={Springer}
}

@article{houtekamer2005ensemble,
  title={Ensemble kalman filtering},
  author={Houtekamer, Peter L and Mitchell, Herschel L},
  journal={Quarterly Journal of the Royal Meteorological Society: A journal of the atmospheric sciences, applied meteorology and physical oceanography},
  volume={131},
  number={613},
  pages={3269--3289},
  year={2005},
  publisher={Wiley Online Library}
}

@article{Leeuwen2020,
author = {Peter Jan van Leeuwen},
title = {A consistent interpretation of the stochastic version of the Ensemble Kalman Filter},
journal = {Quarterly Journal of the Royal Meteorological Society},
volume = {146},
number = {731},
pages = {2815-2825},
doi = {https://doi.org/10.1002/qj.3819},
eprint = {https://rmets.onlinelibrary.wiley.com/doi/pdf/10.1002/qj.3819},
year = {2020}
}

@article{evensen1994sequential,
title={Sequential data assimilation with a nonlinear quasi-geostrophic model using Monte Carlo methods to forecast error statistics},
author={Geir Evensen},
journal={Journal of Geophysical Research: Oceans},
volume={99},
number={C5},
pages={10143--10162},
year={1994},
publisher={Wiley Online Library}
}

@Article{Chatt2020,
AUTHOR = {A. Chattopadhyay and P. Hassanzadeh and D. Subramanian},
TITLE = {Data-driven predictions of a multiscale Lorenz 96 chaotic system using machine-learning methods: reservoir computing, artificial neural network, and long short-term memory network},
JOURNAL = {Nonlinear Processes in Geophysics},
VOLUME = {27},
YEAR = {2020},
NUMBER = {3},
PAGES = {373--389},
URL = {https://npg.copernicus.org/articles/27/373/2020/},
DOI = {10.5194/npg-27-373-2020}
}

@article{LeCun2015,
    author = {Y. LeCun and Y. Bengio and G. Hinton},
    title = {Deep Learning},
    journal = {Nature},
    number = {521},
    year = {2015},
    pages = {436}
}

@article{Chai2021,
title = {Deep learning in computer vision: A critical review of emerging techniques and application scenarios},
journal = {Machine Learning with Applications},
volume = {6},
pages = {100134},
year = {2021},
issn = {2666-8270},
doi = {https://doi.org/10.1016/j.mlwa.2021.100134},
author = {J. Chai and H. Zeng and A. Li and E. W.T. Ngai}
}

@INPROCEEDINGS{Radosavovic2020,
  author={I. Radosavovic and R. P. Kosaraju and R. Girshick and K. He and P. Doll\'{a}r},
  booktitle={2020 IEEE/CVF Conference on Computer Vision and Pattern Recognition (CVPR)}, 
  title={Designing Network Design Spaces}, 
  year={2020},
  pages={10425-10433},
  doi={10.1109/CVPR42600.2020.01044}
}

@ARTICLE{Ye2020,
  author={L. Ye and Z. Liu and Y. Wang},
  journal={IEEE Transactions on Multimedia}, 
  title={Dual Convolutional LSTM Network for Referring Image Segmentation}, 
  year={2020},
  volume={22},
  number={12},
  pages={3224-3235},
  doi={10.1109/TMM.2020.2971171}
}

@inproceedings{Tan2019,
  title={Efficientnet: Rethinking model scaling for convolutional neural networks},
  author={M. Tan and Q. Le},
  booktitle={International conference on machine learning},
  pages={6105-6114},
  year={2019},
  organization={PMLR}
}

@INPROCEEDINGS{Chai2019,
  author={J. Chai and A. Li},
  booktitle={2019 International Conference on Machine Learning and Cybernetics (ICMLC)}, 
  title={Deep Learning in Natural Language Processing: A State-of-the-Art Survey}, 
  year={2019},
  pages={1-6},
  doi={10.1109/ICMLC48188.2019.8949185}
}

@ARTICLE{Otter2021,
  author={D. W. Otter and J. R. Medina and J. K. Kalita},
  journal={IEEE Transactions on Neural Networks and Learning Systems}, 
  title={A Survey of the Usages of Deep Learning for Natural Language Processing}, 
  year={2021},
  volume={32},
  number={2},
  pages={604-624},
  doi={10.1109/TNNLS.2020.2979670}
}

@article{Benes2017,
author = {K. Bene\v{s} and M. K. Baskar and L. Burget},
title = {Residual Memory Networks in Language Modeling: Improving the Reputation of Feed-Forward Networks},
journal = {Proc. Interspeech},
pages = {284-288},
year = {2017},
doi = {10.21437/Interspeech.2017-1442}
}

@article {Hakim2024,
author = {G. J. Hakim and S. Masanam},
title = {Dynamical Tests of a Deep Learning Weather Prediction Model},
journal = {Artificial Intelligence for the Earth Systems},
year = {2024},
publisher = {American Meteorological Society},
address = {Boston MA, USA},
volume = {3},
number = {3},
doi = {10.1175/AIES-D-23-0090.1},
pages= {e230090},
url ={https://journals.ametsoc.org/view/journals/aies/3/3/AIES-D-23-0090.1.xml}
}

@article{Beucler2021,
   author    =  {T. Beucler and M. Pritchard and S. Rasp and J. Ott and P. Baldi and P. Gentine},
   title     =  {Enforcing analytic constraints in neural networks emulating physical systems},
   year      =  {2021},
   journal   =  {Phys. Rev. Lett.},
   volume    =  {126},
   pages     =  {098302}
}

@inproceedings{song2019generative,
author={Y. Song and S. Ermon},
title={Generative Modeling by Estimating Gradients of the Data Distribution},
booktitle={Proc. NIPS 2019},
pages={11918--11930},
year={2019},
month={Dec.},
address={Vancouver, Canada}
}

@article{Jiao2025,
author = {A. Jiao and H. He and R. Ranade and J. Pathak and Lu Lu},
title = {One-shot learning for solution operators of partial differential equations},
journal = {Nat. Commun.},
volume = {16},
pages = {8386},
year = {2025}
}

@article{Chung2014,
author = {J. Chung and C. Gulcehre and K. Cho and Y. Yere},
year = {2014},
title = {Empirical Evaluation of Gated Recurrent Neural Networks on Sequence Modeling},
DOI = {arXiv:1412.3555.}
}

@inproceedings{Hao2024,
 author = {Z. Hao and J. Yao and C. Su and H. Su and Z. Wang and F. Lu and Z. Xia and Y. Zhang and S. Liu and Lu Lu and J. Zhu},
 booktitle = {Advances in Neural Information Processing Systems},
 doi = {10.52202/079017-2442},
 editor = {A. Globerson and L. Mackey and D. Belgrave and A. Fan and U. Paquet and J. Tomczak and C. Zhang},
 pages = {76721-76774},
 publisher = {Curran Associates, Inc.},
 title = {PINNacle: A Comprehensive Benchmark of Physics-Informed Neural Networks for Solving PDEs},
 volume = {37},
 year = {2024}
}

@article{JIANG2024,
title = {Fourier-MIONet: Fourier-enhanced multiple-input neural operators for multiphase modeling of geological carbon sequestration},
journal = {Reliability Engineering \& System Safety},
volume = {251},
pages = {110392},
year = {2024},
issn = {0951-8320},
doi = {https://doi.org/10.1016/j.ress.2024.110392},
author = {J. Zhongyi and M. Zhu and Lu Lu},
}

@article{ZHU2023,
title = {Fourier-DeepONet: Fourier-enhanced deep operator networks for full waveform inversion with improved accuracy, generalizability, and robustness},
journal = {Computer Methods in Applied Mechanics and Engineering},
volume = {416},
pages = {116300},
year = {2023},
issn = {0045-7825},
doi = {https://doi.org/10.1016/j.cma.2023.116300},
author = {M. Zhu and S. Feng and Y. Lin and Lu Lu}
}

@article{Lu2021b,
author = {Lu Lu and X. Meng  and Z. Mao and G. E. Karniadakis},
title = {DeepXDE: A Deep Learning Library for Solving Differential Equations},
journal = {SIAM Review},
volume = {63},
number = {1},
pages = {208-228},
year = {2021},
doi = {10.1137/19M1274067},
eprint = { 
https://doi.org/10.1137/19M1274067}
}

@article{lipton2015,
  title={A critical review of recurrent neural networks for sequence learning},
  author={Z. C. Lipton and J. Berkowitz and C. Elkan},
  journal={arXiv preprint arXiv:1506.00019},
  year={2015}
}

@book{Graves2012Supervised,
  author    = {Alex Graves},
  title     = {Supervised Sequence Labelling with Recurrent Neural Networks},
  publisher = {Springer},
  address   = {Berlin, Heidelberg},
  series    = {Studies in Computational Intelligence},
  year      = {2012},
  doi       = {10.1007/978-3-642-24797-2}
}

@inproceedings{Pascanu2013,
  title={On the difficulty of training recurrent neural networks},
  author={R. Pascanu and T. Mikolov and Y. Bengio},
  booktitle={International conference on machine learning},
  pages={1310-1318},
  year={2013},
  organization={Pmlr}
}

@INBOOK{Kolen2001,
  author={J, F. Kolen and S. C. Kremer},
  booktitle={A Field Guide to Dynamical Recurrent Networks}, 
  title={Gradient Flow in Recurrent Nets: The Difficulty of Learning LongTerm Dependencies}, 
  year={2001},
  publisher = {IEEE},
  pages={237-243},
  doi={10.1109/9780470544037.ch14}}

@article{Graves2012,
  title={Long short-term memory},
  author={Alex Graves},
  journal={Supervised sequence labelling with recurrent neural networks},
  pages={37-45},
  year={2012},
  publisher={Springer}
}

@article{Graves2013,
author = {Alex Graves},
year = {2013},
title = {Generating Sequences With Recurrent Neural Networks},
url = {https://arxiv.org/abs/1308.0850v5}
}

@ARTICLE{Hochreiter1997,
  author={S. Hochreiter and J. Schmidhuber},
  journal={Neural Computation}, 
  title={Long Short-Term Memory}, 
  year={1997},
  volume={9},
  number={8},
  pages={1735-1780},
  url={10.1162/neco.1997.9.8.1735}}

@article{Gers2000,
    author = {F. A. Gers and J. Schmidhuber and F. Cummins},
    title = {Learning to forget: continual prediction with LSTM},
    journal = {Neural computation},
    volume = {12},
    number = {8},
    pages = {2451-2471},
    year =  {2000},
    url = {https://doi.org/10.1162/089976600300015015}
}

@article{LETKF_2007,
  title={Efficient data assimilation for spatiotemporal chaos: A local ensemble transform Kalman filter},
  author={B. R. Hunt and E. J. Kostelich and I. Szunyogh},
  journal={Physica D: Nonlinear Phenomena},
  volume={230},
  number={1-2},
  pages={112--126},
  year={2007},
  publisher={Elsevier},options = {skipbib=true}
}

@book{evensen2009data,
  title={Data assimilation: the ensemble Kalman filter},
  author={Gier Evensen},
  year={2009},
  publisher={Springer Science \& Business Media}
}

@ARTICLE{LETKF_2020,
    AUTHOR = {M. Buehner},
    YEAR = {2020},
    TITLE = {Local Ensemble Transform Kalman Filter with Cross Validation},
    JOURNAL = {MWR},
    VOLUME = {148},
    PAGES = {2265--2282},
    DOI = {10.1175/MWR-D-19-0402.1}
}

@article{Andrieu2010,
author = "C. Andrieu and A. Doucet and R. Holenstein",
title = "Particle Markov chain Monte Carlo methods",
year = "2010",
journal = "J. R. Statist. Soc. B",
volume = "72",
number = "3",
pages = "269-342"
}

@article{Kang2018,
AUTHOR = {K. Kang and V. Maroulas and I. Schizas and F. Bao},
 TITLE = {Improved distributed particle filters for tracking in a wireless sensor network},
JOURNAL = {Comput. Statist. Data Anal.},
FJOURNAL = {Computational Statistics \& Data Analysis},
VOLUME = {117},
  YEAR = {2018},
 PAGES = {90-108},
  ISSN = {0167-9473},
MRCLASS = {Expansion},
MRNUMBER = {3711483}
}

@article{MTAC2012,
title = "A random map implementation of implicit filters",
journal = "J.  Comput. Phys.",
Fjournal = "Journal of Computational Physics",
volume = "231",
number = "4",
pages = "2049-2066",
year = "2012",
 AUTHOR = "M. Morzfeld and X. Tu and E. Atkins and A. J. Chorin",
}

@article{pitt1999filtering,
  title={Filtering via simulation: Auxiliary particle filters},
  author={M. K. Pitt and N. Shephard},
  journal={Journal of the American statistical association},
  volume={94},
  number={446},
  pages={590-599},
  year={1999},
  publisher={Taylor \& Francis},options = {skipbib=true}
}

@article{Gordon1993, 
AUTHOR = {N. J. Gordon and D. J. Salmond and A. F. M. Smith}, 
TITLE = {Novel approach to nonlinear/non-Gaussian Bayesian state estimation},
JOURNAL = {IEE PROCEEDING F},
 VOLUME = {140}, 
 YEAR = {1993}, 
 NUMBER = {2}, 
 PAGES = {107-113}
}

@article{hu_vanLeeuwen_2021,
 AUTHOR = {C.-C. Hu and P. J. van Leeuwen},
YEAR = {2021},
TITLE = {A particle flow filter for high-dimensional system applications},
JOURNAL = {QJRMS},
VOLUME = {147},
PAGES = {2352-2374},
DOI = {10.1002/qj.4028}
}

@ARTICLE{Sny,
AUTHOR="C. Snyder and T. Bengtsson and P. Bickel and J. Anderson",
TITLE="Obstacles to high-dimensional particle filtering",
JOURNAL=" Mon. Wea. Rev.",
 YEAR=2008,
VOLUME=136,
PAGES="4629-4640"
}

@article{Bao2019a,
  author  =      {R. Archibald and F. Bao and X. Tu},
  title   =      {A direct method for parameter estimations},
  journal =      {J. Comput. Phys.},
  volume  =      {398},
  pages   =      {108871},
  year    =      {2019}
}

@article{YLiu2024,
	author={Y. Liu and M. Yang and Z. Zhang and Y. Cao and G. Zhang and F. Bao},
	title ={Diffusion-Model-Assisted Supervised Learning of Generative Models for Density Estimation},
	journal = {Journal of Machine Learning for Modeling and Computing},
	volume = {5},
	number = {1},
	year = {2024},
	pages = {25-38}
}

@inproceedings{YSong2021,
title={Score-Based Generative Modeling through Stochastic Differential Equations},
author={Y. Song and J. Sohl-Dickstein and D. P Kingma and A. Kumar and S. Ermon and B. Poole},
booktitle={International Conference on Learning Representations},
year={2021},
url={https://openreview.net/forum?id=PxTIG12RRHS}
}
\nocite{*}

\end{document}